\documentclass[11pt,a4]{article}

\usepackage{algorithm,algpseudocode}

\usepackage{amsmath,amsthm,amssymb}
\usepackage{geometry,subfig,epsfig,epstopdf,bbm}

\newtheorem{theorem}{Theorem}[section]

\newtheorem{remark}[theorem]{Remark}

\numberwithin{equation}{section}
\newtheorem{exmp}[theorem]{Example}
\usepackage{enumerate}
\usepackage[shortlabels]{enumitem}

\newcommand{\R}{\mathbb{R}}
\newcommand{\bxi}{\boldsymbol{\xi}}

\newcommand{\E}{\mathbb{E}}
\newcommand{\Prm}{\mathbb{P}}
\newcommand{\bta}{\boldsymbol{\eta}}
\newcommand{\balp}{\boldsymbol{\alpha}}
\newcommand{\bgam}{\boldsymbol{\gamma}}
\newcommand{\bbeta}{\boldsymbol{\beta}}

\title{A dynamical polynomial chaos approach for long-time evolution of SPDEs}

\author{H. Cagan Ozen\thanks{Department of Applied Physics \& Applied Mathematics, Columbia University, New York, NY 10027, USA (\href{mailto:hco2104@columbia.edu}{hco2104@columbia.edu}, \href{gb2030@columbia.edu}{gb2030@columbia.edu}). }  \and Guillaume Bal\footnotemark[1] }

\date{\today}

\begin{document}

\maketitle

\begin{abstract}
We propose a Dynamical generalized Polynomial Chaos (DgPC) method to solve time-dependent stochastic partial differential equations (SPDEs) with white noise forcing. The long-time simulation of SPDE solutions by Polynomial Chaos (PC) methods is notoriously difficult as the dimension of the stochastic variables increases linearly with time. Exploiting the Markovian property of white noise, DgPC \cite{OB16} implements a restart procedure that allows us to expand solutions at future times in terms of orthogonal polynomials of the measure describing the solution at a given time and the future white noise. The dimension of the representation is kept minimal by application of a Karhunen--Loeve (KL) expansion. Using frequent restarts and low degree polynomials on sparse multi-index sets, the method allows us to perform long time simulations, including the calculation of invariant measures for systems which possess one.  We apply the method to the numerical simulation of stochastic Burgers and Navier Stokes equations with white noise forcing. Our method also allows us to incorporate time-independent random coefficients such as a random viscosity.  We propose several numerical simulations and show that the algorithm compares favorably with standard Monte Carlo methods.

\end{abstract}

\pagestyle{myheadings}
\thispagestyle{plain}
\markboth{ H.C. Ozen and G. Bal}{DgPC}


\section{Introduction}
\label{}

Stochastic partial differential equations (SPDEs), or partial differential equations (PDEs) with uncertainties, play an important role in many areas of engineering and applied sciences such as turbulence, flows in random media, solid mechanics, filtering, and finance. Numerical simulations of SPDEs are typically based on the Monte Carlo (MC) or the  Polynomial Chaos (PC) methods  \cite{JK09,LPS14, GS91, LeMK10}. In both cases, the long-time simulation of SPDEs proves to be quite expensive~\cite{HLRZ06,XK03,LKHG01}. In this paper, we extend the PC-based method developed for stochastic differential equations in \cite{OB16} to the setting of SPDEs.

The PC method originates from the works \cite{W38,CM47} and enables us to expand (square integrable) functionals of Brownian motion in a basis of Hermite polynomials. Applications of such a framework to random flows and turbulence theory are examined in  \cite{C74,OB67}. More recently, works in \cite{GS91,GRH99} combined the PC method with the Karhunen--Loeve (KL) expansion  \cite{K47,L77}  to study structural mechanics problems.  The generalized polynomial chaos (gPC) developed in \cite{XK02, XK03} next extends the Hermite PC to a set of non-Gaussian random parameters. 

The main advantage of the PC method is that it allows us to propagate stochasticity by providing expansions of quantities of interest in terms of appropriate uncertainties while in effect replacing the stochastic equations by systems of deterministic equations. Once these deterministic equations are solved, the statistical properties of the solution can be readily inferred from the coefficients of the expansion, which facilitates uncertainty quantification. In some cases, the PC method can propagate uncertainties with a substantially lower cost than MC methods; especially for low dimensional uncertainties~\cite{GS91, XK02,XK03,HLRZ06,BM13, OB16}; see also \cite{BalUQ16}. However, in cases of high dimensional random parameters, the efficiency of the PC method is reduced because of the large number of terms that appear in the expansion. The presence of a temporal random forcing is a serious challenge to the PC method as the number of stochastic variables increases linearly with time, which hinders the capability of the PC method for long-time computations \cite{OB16,HLRZ06}. Moreover, the optimality of the initial basis typically degrades as time evolves due to the presence of nonlinearities in the equation. These drawbacks were addressed in \cite{WK05,WK06, HLRZ06,Luo,GSVK10, PL2012,BM13}. In our previous work \cite{OB16}, we proposed the \textit{Dynamical generalized Polynomial Chaos} (DgPC) method to solve low dimensional stochastic differential equations (SDEs) driven by white noise. The novelty of the DgPC algorithm is that it utilizes a restart procedure and constructs polynomial chaos expansions (PCEs) dynamically in time by using orthogonal polynomials of the projections of the solution at the current steps and the future random forcing~\cite{OB16}. This allows the algorithm to keep reasonably sparse random bases, and to mitigate the aforementioned curse of dimensionality and loss of optimality. Numerical experiments in that reference illustrate that such a restart procedure has the ability to accurately estimate long-time solutions of SDEs; see the relevant computational and theoretical details in \cite{OB16}.

Although the algorithm performs well for low dimensional SDEs, extension to larger systems is challenging and requires serious modifications. In this manuscript, we extend DgPC  to the framework of SPDEs driven by white noise. While solutions to SPDEs are, in general, high dimensional random fields, they may lend themselves in some cases to low dimensional representations \cite{GS91,BPL93, JSK02,ST06,HLRZ06,DGR07, LeMK10}. Armed with this fact, we propose at each restart to use the KL expansion to compress the solution into a representation involving a finite number of random modes, i.e., a projection onto a lower dimensional manifold. In cases where the modeling equations contain non-forcing random inputs other than Brownian forcing, such as a random viscosity, the KL expansion is applied to the solution and the random parameters together so that the algorithm automatically selects the intrinsic variables, which have the largest influence on the solution. 

Since KL expansions are known to be optimal in the mean square sense, at each restart point in time, only a few dominating, most energetic, random modes are chosen and incorporated into PCE to represent the future solution. However, the combination of the random KL modes and the random forcing variables brings about high dimensionality. The computational challenges then become: (i) computing orthogonal polynomials of arbitrary multivariate distributions; (ii) keeping the number of terms in the expansion as small as possible. 

The construction of orthogonal polynomials of evolving multivariate distributions is possible by estimating their statistical moments~\cite{gautschi,OB16}, which is, in general, a computationally intensive procedure. In this work, we estimate the moments of the solution using its PCE through a sampling approach to greatly reduce their computational cost; see \cite{AGPRH12} for a similar sampling methodology. We also show that the method is robust with respect to re-sampling.

The second challenge (ii) is a major problem for all PC-based methods. 
The KL decomposition is computationally expensive as it requires solving a large eigenvalue problem. For problems of moderate size, we find that the eigenvalue problem is solved efficiently by using a Krylov subspace method. For larger problems, in order to mitigate both memory requirements and computational cost of the KL expansion, we find low-rank approximations to the covariance matrices based on their PC representations and without assembling them. The algorithm leverages randomized projection methods as introduced in \cite{MRM11,HMT11}, and uses appropriate random matrices to obtain fast and accurate solutions of large eigenvalue problems. 
After selecting the dominating modes in the KL expansion, we make use of the sparse truncation technique proposed in \cite{HLRZ06,Luo} to further reduce the number of terms in a PCE. For long-time computations,  we also develop an adaptive restart scheme, which adapts the time lag between restarts based on the acuity of the nonlinear effects. 

The use of compression techniques to exploit intrinsic lower dimensional structures of solutions of SPDEs is not new and is in fact necessary in many contexts; see \cite{BPL93,DGR07, VWK08, SL09, CSK13,CHZ13}. The novelty of our approach is that a lower dimensional representation of the solution is learned online and integrated into a PCE to integrate future random forcing and represent future solutions. This procedure is computationally viable and the combination of the aforementioned  ideas allows us to attain a reasonable accuracy in the  long-time evolutions of SPDEs for a reasonable computational cost.

As we are interested in the long-time evolution of SPDEs, we restrict ourselves here to dynamics with a dissipation mechanism. Equilibrium statistics and asymptotic properties of the solutions are relevant in many applications and have been extensively studied in the literature \cite{OB67, C74, S91, HLRZ06,HM06, BM13}. Based on these motivations, we provide numerical experiments for a 1D randomly forced Burgers equation and a 2D stochastic Navier--Stokes (SNS) system. All equations are driven by white noise and satisfy periodic (spatial) boundary conditions. To demonstrate the efficiency of our algorithm, we present both short- and long-time computations. In some cases, we model the viscosity as a random process.  Statistical moments obtained by the algorithm are compared to the standard MC methods for short times to assess the accuracy of the algorithm. Results show promising speed-ups over standard MC methods. We exhibit convergence behavior in terms of the degree of the expansion, the number of random modes retained in the KL expansion, and the (adaptive) restart time. In some cases, we provide a purely PCE-based numerical verification of the convergence of the process to its invariant measure.

The outline of the paper is as follows. Section \ref{sec:KL_PC} introduces the basic setting for the PC and KL expansions. In the following section, the main algorithm and its implementation details are discussed. Numerical experiments are conducted in section \ref{sec:Numerical} using several settings for the stochastic Burgers equation and Navier--Stokes system. Some concluding remarks are offered in section \ref{sec:summary}.

\section{An overview of Karhunen--Loeve  and polynomial chaos expansions} \label{sec:KL_PC}

Throughout the manuscript, we consider the following time-dependent stochastic partial differential equation (SPDE) driven by the white noise $\dot{W}(t,\omega)$: 
\begin{align}\label{eq:gen_spde}
\left\{
\begin{aligned}
& \partial_t  u(x,t,\omega) = \mathcal{L}(u(x,t,\omega)) + \sigma(x) \, \dot{W}(t,\omega), \quad x \in G \subset \R^d, \, \omega \in \Omega, \, t \in [0,T], \\
&u(x,0,\omega)= u_0(x,\omega), \quad x \in G, \, \omega \in \Omega, 
\end{aligned}
\right.
\end{align}
where, for concreteness, $G$ is the $d$-dimensional torus so that $u$ and its derivatives are periodic functions in the variable $x$. The DgPC algorithm may easily be extended to more general boundary conditions and geometries. 
Above,  $\mathcal{L}$ is a possibly non-linear differential operator in the spatial variables.
The solution takes values in $\R^p$. In the numerical simulations, the parameters $d$ and $p$ are set to either $1$ or $2$.

\subsection{Karhunen--Loeve expansion}

Given a probability space $(\Omega, \mathcal{F}, \Prm)$, where $\Omega$ is a sample space equipped with the sigma-algebra $\mathcal{F}$ and the probability measure $\Prm$, we denote by $L^2(G \times \Omega)$ the Hilbert space of square integrable random fields on $G$. For a random field $u \in L^2(G \times \Omega)$, we define the expectation 
\[
\bar{u}(x) = \E [u(x,\omega)] = \int_{\Omega} u(x,\omega) \Prm(d\omega),
\]
and the covariance
\begin{align} \label{eq:eig_C}
\mbox{Cov}_u(x,y) = \E [ (u(x,\omega) - \bar{u}(x) )  \,(u(y,\omega) - \bar{u}(y))' ], \quad x,y \in G,
\end{align}
where $'$ denotes the transpose.

The KL expansion of a time-independent random field $u \in L^2(G \times \Omega)$ with a continuous covariance is as follows:  
\begin{align} \label{eq:KL}
u(x,\omega) = \bar{u}(x) + \sum_{l=1}^{\infty} \sqrt{\lambda_l} \,  \eta_l(\omega) \,  \phi_l(x),
\end{align}
where the eigenvalues $\lambda_l$ are nonnegative, and the eigenfunctions $\phi_l$ and the mean zero random variables $\eta_l$ are given by
\begin{align*}
\langle \mbox{Cov}_u(x, \cdot) \, ,\phi_l \rangle_{L^2(G)} = \lambda_l \, \phi_l(x),  \qquad   \eta_l(\omega) = \lambda_l^{-1/2} \langle [u(\cdot, \omega)- \bar{u}] \, ,\phi_l \rangle_{L^2(G)}.
\end{align*}
Here, $\langle \cdot, \cdot \rangle_{L^2(G)}$ denotes the spatial inner product on $G$. The spatial and random modes satisfy bi-orthogonality, i.e. $\E[\eta_l \, \eta_k] = \delta_{lk}$ and $\langle \phi_l , \phi_k \rangle_{L^2(G)}=\delta_{lk}$, where $\delta_{lk}$ denotes the Kronecker delta function;  see \cite{K47,L77}. 

The major feature of the KL decomposition is that the truncation after a finite number, denoted by $D$ hereafter, of terms is optimal in the $L^2$ sense. How $D$ should be chosen obviously depends on the spectrum of the covariance operator. When the process shows a high degree of correlation, then typically $D$ may be chosen relatively small due to the rapid decay of the eigenvalues  \cite{GS91, LeMK10, ST06}. This makes the KL expansion a useful dimensionality reduction technique in many applications and will play a crucial role in our algorithm to compress the dimensionality of stochastic solutions of \eqref{eq:gen_spde}.  

\subsection{Polynomial chaos expansion} 

We assume that the random solution $u(x,t,\omega)$ of \eqref{eq:gen_spde} is a second-order stochastic process for all $t \in [0,T]$, and for simplicity, that the initial condition is deterministic. The PC method for \eqref{eq:gen_spde} is constructed as follows. 

Given a complete countable orthonormal system $m_i(t) \in L^2[0,T]$, $i \in \mathbb{N}$, we project the Brownian motion $W(t,\omega)$ onto $L^2 [0,T]$ by defining $\xi_i(\omega) := \int_0^T m_i(t)\, dW(t,\omega) $. Then,  $\boldsymbol{\xi}= (\xi_1,\xi_2,\ldots)$ consists of a countable number of  independent and identically distributed standard Gaussian random variables, and the convergence 
\begin{align} \label{eq:conv_W}
 \E \left[ W(t) - \sum_{i=1}^K  \xi_i \, \int_0^t m_i(\tau) \, d\tau \right]^2  \rightarrow 0 , \quad K \rightarrow \infty,
\end{align}
holds for all $t \leq T$. Basic examples of orthonormal systems are trigonometric functions and wavelets; see ~\cite{HLRZ06, Luo,MNGK04}. 

The SPDE solution is a nonlinear functional of Brownian motion on the interval $[0,T]$. Therefore, the projection onto $\bxi$ enables us to consider the solution depending on a countable number of random variables, i.e., $u(x,t) \in L^2(\Omega,\sigma(\boldsymbol{\xi}) , \Prm)$. Then, the polynomial chaos expansion (PCE) of the solution $u(x,t)$ at $t \in [0,T]$ is given by
\begin{align} \label{eq:PCE}
u(x,t,\boldsymbol{\xi}) = \sum_{\balp \in \mathcal{J}} u_{\balp}(x,t)\,  T_{\balp}(\boldsymbol{\xi}), \qquad  u_{\balp}(x,t) = \E[u(x,t,\boldsymbol{\xi}) \, T_{\balp}(\boldsymbol{\xi})]. 
\end{align}
 Note that the expansion is a sum of orthogonal projections onto orthogonal subspaces spanned by the Wick polynomials  $T_{\balp}(\bxi) := \prod_{i=1}^{\infty} H_{\alpha_i}(\xi_i)$, where $H_{n}$ is the $n$th degree normalized one-dimensional Hermite polynomial and $\balp$ belongs to set of multi-indices with finite number of nonzero components
$
\mathcal{J} = \{ \balp=(\alpha_1,\alpha_2, \ldots) \, | \, \alpha_j \in \mathbb{N} \cup \{0\}, \, |\balp|=\sum_{i=1}^{\infty} \alpha_i < \infty \}.
$
The completeness of the orthonormal polynomials $T_{\balp}(\bxi)$ in $L^2$, i.e. convergence of the expansion \eqref{eq:PCE}, is established in~\cite{CM47} and the result is known as the Cameron-Martin theorem. In the following, we will use the term Hermite PCE to emphasize that the expansion utilizes Gaussian random variables. 

All statistical information of the random field is contained in the deterministic coefficients $u_{\balp}$ of the expansion \eqref{eq:PCE}. Specifically, the first two moments are given by
\begin{align} \label{eq:first_mom}
\E [u] = u_{\boldsymbol{0}}(x,t), \qquad \E[u^2] = \sum_{\balp \in \mathcal{J}} u_{\alpha}^2(x,t).
\end{align}
Higher order moments can be obtained using the triple products $\E [T_{\balp}(\bxi) T_{\bgam}(\bxi) T_{\bbeta}(\bxi)]$ or the Hermite product formula for $\bxi$; see \cite{MR04,HLRZ06, Luo, OB16}. 

In practice, we truncate the expansion \eqref{eq:PCE} using polynomials of maximum degree $N$ and $K$ number of random variables $(\xi_1,\ldots,\xi_K)$, which results in the following approximation
\begin{align} \label{eq:simp_trunc}
u \approx u_{K,N}(x,t,\xi_1, \ldots ,\xi_K) := \sum_{|\balp| \leq N} u_{\balp}(x,t) \prod_{i=1}^K H_{\alpha_i}({\xi_i}).
\end{align}
This simple truncation gives rise to $K+N \choose{K}$ terms in the approximation. 

Two major problems of the PCE in this context are as follows. To maintain a prescribed accuracy as time evolves, the degree of freedom $K$ should grow in accordance with the error behavior of \eqref{eq:conv_W}, which brings about high dimensional randomness that needs to be incorporated in the expansion. Therefore, the computational cost increases rapidly with dimension, which decreases the efficiency of the PC method. The second related problem of PC is that initial predetermined bases may lose their optimality and even fail to converge for long-time evolutions~\cite{WK05, WK06,HLRZ06,BM13,MNGK04,GSVK10, OB16}.

\subsection{Dynamical generalized PC} \label{sec:DgPC}

To alleviate the aforementioned bottlenecks of the PC method, namely high dimensionality and long-term integration, several approaches have been proposed~\cite{XK02, Xiu_thesis, WK05,WK06, HLRZ06, GSVK10, PL2012, HS14}. In the following, we briefly summarize the approach called \textit{Dynamical generalized Polynomial chaos }(DgPC), which was developed in \cite{OB16} to solve (reasonably low dimensional) stochastic differential equations (SDEs) driven by Brownian motion over long times.

The solution $u(t,\omega)$ of a $d$-dimensional SDE system can be seen as a nonlinear functional of the initial condition $u_0$ and the countable number of variables $\bxi$. The approach in \cite{OB16} hinges upon two important observations: (i) for sufficiently small time lags $\varepsilon>0$, the solution $u(t+\varepsilon)$, can be efficiently captured by low-order polynomials in $u(t)$; and (ii) the solution $u(t+\varepsilon)$  depends on Brownian motion only on the interval $s \in [t,t+\varepsilon]$ by the Markov property. These two remarks allowed us to employ evolving PCEs based on polynomials of the random forcing and projections of the solution onto underlying chaos spaces. In this way, the algorithm can ``forget" about the past and PC representations remain reasonably sparse in the evolving basis while the accuracy is retained for long times. 

To be more specific, and given a set of increasing restart times $t_j$ on $(0,T)$, $j=1,\ldots,n$ and $n  \in \mathbb{N}$,  the algorithm projects the random variable $u(t_j)$ onto a finite chaos space at each $t_j$. This projection will be denoted by $u_j$ in the following. A major challenge in the implementation is to construct orthogonal polynomials of the varying distribution of the random variable $u_j$. The PCE \eqref{eq:PC_exp} gives rise to natural formulas \eqref{eq:first_mom} to compute statistical moments of the random variable. The algorithm constructs orthogonal polynomials of $u_j$ based on the knowledge of its moments via a Gram--Schmidt procedure. To obtain evolution equations for the expansion coefficients, the algorithm then computes necessary triple products for $u_j$ and performs Galerkin projection onto subspaces spanned by orthogonal polynomials in $u_j$ and $\bxi$. Convergence of the algorithm is related to the unique solvability of the moment problem of the measures associated to variables $u_j$; see \cite{BC81, EMSU12}.   Related theoretical aspects and the estimates of computational costs can be found in detail in \cite{OB16}. 

In \cite{OB16}, the algorithm is applied to low dimensional nonlinear SDE systems and provides accurate long-time simulations with a small cost compared to standard Hermite PCE. In particular, we computed invariant measures of SDEs in some cases, which standard Hermite PCE typically fails to achieve. In the subsequent sections, we propose an extension of the DgPC algorithm to SPDEs, which requires several key modifications.

\section{Dynamical generalized polynomial chaos for SPDEs} \label{sec:DgPC_spde}

In this section, we present the algorithmic details of the DgPC method applied to the general SPDE \eqref{eq:gen_spde}. The main ideas are summarized as follows.

Let a decomposition  $0= t_0 < t_1<\ldots<t_n=T$ be given. Following our discussion in section \ref{sec:DgPC} and using the Markov property, the solution $u(x,t_{j+1},\omega)$, $ 0 \leq j< n$, can be  
represented in a PC expansion in terms of $u(x,t_{j},\omega)$ and $\bxi_j$, where $\bxi_j=(\xi_{j,1},\xi_{j,2},\ldots)$ denotes the Gaussian random variables
required for Brownian forcing on the interval $[t_j , t_{j+1}]$. Let $u_j(x,t_j,\omega)$ denote the projection of the solution $u(x,t_j,\omega)$ onto the polynomial chaos space. We will also use the shorthand notation $u_j$ to denote this projection. To construct a PCE in terms of polynomials of $u_j$,  we separate the spatial dependence and randomness via the KL expansion \eqref{eq:KL}:
\begin{align} \label{eq:KL_PC}
u_j = \bar{u}_j (x) + \sum_{l=1}^{\infty} \sqrt{\lambda_{j,l}} \, \eta_{j,l}(\omega) \, \phi_{j,l}(x).
\end{align}
Let $\bta_j := (\eta_{j,1},\eta_{j,2},\ldots)$ denote the random modes. Since the solution $u(x,t,\omega)$, $t \in [t_j,t_{j+1}]$, is a functional of the random forcing $\bxi_j$ and modes $\bta_j$, the next step PCE  $u_{j+1}(x,t,\omega)$ is given by
\begin{align} \label{eq:PC_exp}
u_{j+1}(x,t,\omega)  = \sum_{\balp \in \mathcal{J}} u_{j+1,\balp}(x,t) \, T_{\balp} (\bxi_j(\omega) , \bta_j(\omega)), \quad t \in [t_j,t_{j+1}],
\end{align}  
with the notation $u_{j+1}(x,t_j,\omega) = u_j$, where $T_{\balp}$ denotes an orthonormal basis in its arguments. The expansion dynamically needs a PC basis depending on the random forcing and evolving random KL modes of the solution. The coefficients $u_{j+1,\balp}(x,t)$ satisfy a PDE system obtained by Galerkin projection of the SPDE \eqref{eq:gen_spde} onto the space spanned by $T_{\balp}(\bxi_j,\bta_j)$. Statistical properties can be retrieved after solving the induced PDE system.  

Here are the computational bottlenecks of this approach: (i) the simple truncation \eqref{eq:simp_trunc} yields a large number of terms in the expansion, which leads to long computational times to solve the deterministic evolution equations; (ii) estimating the terms appearing in the  KL expansion \eqref{eq:KL_PC} is a major computational bottleneck, especially in higher spatial dimensions; and (iii) computation of the orthogonal polynomials $T_{\balp}(\bxi_{j},\bta_j)$ may also require intensive amount of computation. In the following, we address these issues in turn.  
  
\subsection{Sparse truncation}

The number of terms in the simple truncation \eqref{eq:simp_trunc} in the Hermite PCE increases rapidly with respect to $N$ and $K$. Given sufficient regularity of the solution, the expansion coefficients decay both in the number of Gaussian variables $K$ and the degree of polynomials $N$. This observation led the authors \cite{HLRZ06} to introduce a sparse truncation of the multi-index set and retain a truncated random basis, which keeps lower (higher) order polynomials in $\xi_i$ with larger (smaller) subscripts. This truncation can be quantified using an estimate for the decay rate of the coefficients; see \cite{Luo, HLRZ06, XL09, FST05, BS11}. Following \cite{Luo,HLRZ06, FST05}, we introduce a sparse index
\[
r = (r_1,\ldots,r_K), \quad N \geq r_1 \geq r_2 \geq \ldots \geq r_K, 
\]
and define the corresponding sparse multi-index set 
\begin{align} \label{eq:sparse_mi}
J_{K,N}^r := \{ \balp=(\alpha_1,\ldots,\alpha_K), \, |\balp|\leq N , \alpha_i \leq r_i \}.
\end{align}
Basically, the index $r$ keeps track of how much degree we want in each variable $\xi_i$. Using  \eqref{eq:sparse_mi}, one can define the corresponding version of the PC expansion \eqref{eq:PCE} which might have drastically reduced number of terms. This is possible by a suitable choice of the index $r$ so that ineffective cross terms in high degree polynomials are eliminated. 

\subsection{KL expansion} \label{sec:KLexpansion}

At each restart $t_j$, we employ the KL expansion \eqref{eq:KL_PC} for the projected random field $u_j$, and the expansion is truncated after a finite number of $D$ terms. The decomposition yields the eigenvalues ${\lambda_{j,l}}$ and the eigenfunctions $\bta_j = (\eta_{j,l})$ for $l=1,\ldots,D$. Therefore, in addition to the Gaussian variables $\bxi_j$, we have the random modes $\bta_j$ in the PCE \eqref{eq:PC_exp} so that the total number of random variables becomes $K+D$. To accommodate $\bta_j$, we extend the multi-index set \eqref{eq:sparse_mi} to $\mathcal{J}^r_{K+D,N}$. This can also be done by the tensorization $\mathcal{J}^r_{K,N} \otimes \mathcal{J}^r_{D,N}$ of the multi-index sets since $\bxi_j$ and $\bta_j$ are independent. However, since tensorization yields higher number of terms in the PCE for most values of $K$ and $D$, it is not considered in the following.  

Assuming that the orthonormal basis $T_{\balp}(\bxi_{j-1},\bta_{j-1})$ is constructed in the previous step, the solution $u(x,t_j,\omega)$ is approximated using the  truncated PCE 
\begin{align} \label{eq:PCtruncation}
 u_j = \sum_{\balp \in \mathcal{J}^r_{K+D,N}} u_{j,\balp}(x) \, T_{\balp}(\bxi_{j-1}, \bta_{j-1}),
\end{align}
where the time dependence of the coefficients  $u_{j,\balp}(x)$ is omitted for brevity. To avoid confusion, we note that both the projection \eqref{eq:PC_exp} and its truncation \eqref{eq:PCtruncation} will be denoted by $u_j$. Armed with this approximation, using the orthogonality of random bases, the covariance of $u_j$ is easily estimated  by
\begin{align} \label{eq:cov_pce}
\mbox{Cov}_{u_j}(x,y) = \sum_{\balp>0}  u_{j,\balp}(x) \,  u_{j,\balp}(y)', \quad x,y \in G, 
\end{align} 
where $\balp \in \mathcal{J}^r_{K+D,N}$. 

In practice, we consider a discretization of the spatial domain $G$ with an even mesh parameter $M \in \mathbb{N}$ and grid points $x_m$, $m=1,\ldots,M^d$. Denote by $C$ the resulting covariance matrix. In general, we can use a spectral method, e.g., Fourier series in the case of periodic functions, to approximate the coefficients $u_{j,\balp}$ as
\begin{align} \label{eq:cov_spectral}
u_{j,\balp} (x) \approx \sum_{k=1}^{M^d}  \hat{u}_{j,{\balp}}(k) \, \varphi_k(x),
\end{align}
where $\varphi_k$ are orthogonal global basis functions on $G$ and $\hat{u}_{j,\balp}(k) = \langle u_{j,\balp}, \varphi_k  \rangle_{L^2(G)}$; see \cite{HGG07}. Thus, $u_{j,\balp}(x)$ is approximated by a vector $(u_{j,\balp}(x_m))_{m=1}^{M^d}$ on the grid. Therefore, the dimension of the covariance matrix $C$ becomes of order $\mathcal{O}(M^d \times M^d)$. 

After computing the covariance matrix, the corresponding eigenvalue problem can be solved for the first $D$ largest eigenvalues. Then, the random modes $\bta_j$ are given as
\begin{align} \label{eq:etaij}
\eta_{j,l} = \frac{1}{\sqrt{\lambda_{j,l}}}  \langle \, (u_j-\E[u_j]) \,, \phi_{j,l} \rangle_{L^2(G)} = \frac{1}{\sqrt{\lambda_{j,l}}} \sum_{\balp>0}  \langle u_{j,\balp}, \phi_{j,l} \rangle_{L^2} \, T_{\balp}(\bxi_{j-1}, \bta_{j-1}), \quad l=1,\ldots,D. 
\end{align}
This representation yields the random modes $\bta_j$ as a function of $\bxi_{j-1}$ and the previous modes $\bta_{j-1}$. Here, we assume that the integrals $ \langle u_{j,\balp}, \phi_{j,l} \rangle_{L^2}$
can be computed by an accurate quadrature method. 

\begin{remark} \rm
When the solution $u$ is more than one-dimensional, several implementations of the KL expansion can be considered. For instance, we may apply the KL expansion to each component of the solution $u$ separately and incorporate the resulting individual random modes into one PCE. Although this approach certainly makes the KL step faster, we found that it yielded inaccurate results in DgPC and needed a large number of variables in the PCE to represent the solution accurately since cross covariance structures between the components of the solution are lost. Therefore, in the following, we implement the KL expansion directly to the multi-dimensional solution and produce one set of random modes $\bta$ which represents all of its components. 
\end{remark}

Depending on the resolution of the discretization of the domain $G$ and the dimension $d$, assembling the covariance matrix and solving the corresponding eigenvalue problem may prove dauntingly expensive.  Several methods have been devised to reduce computational costs, such as fast Fourier techniques (e.g., in the case of stationary covariance kernels) or sparse matrix approximations together with Krylov subspace solvers \cite{ST06, GL96,EEU07, KLM09}. 

Here, the assembly of the covariance matrix is performed at each restart via the summation formula \eqref{eq:cov_pce}. In our one dimensional simulations, with $d=1$, this assembly can be carried out reasonably fast. Since the solution of the eigenvalue problem is required only for the  number $D\ll M^d$ of eigenvalues and eigenfunctions, Krylov subspace methods \cite{S11} perform well. We also utilize the implicitly restarted Arnoldi method to efficiently find the few largest eigenvalues and corresponding eigenfunctions; \cite{A51,LSY98}. 

In the higher dimensional simulations, when $d>1$, computing and storing such large covariance matrices become challenging. Covariance matrices are, in general, not sparse and require $\mathcal{O}(M^{2d})$ units of storage in the memory. Moreover, assembling a large matrix at every restart is computationally very expensive for long-time simulations. The problem of computing and storing a large covariance matrix resulting from the KL expansion is not new and was addressed before in \cite{ST06,EEU07,KLM09,CHZ13}. It was noted that although the covariance matrices were dense, they were usually of low-rank; \cite{ST06,EEU07}. Based on this observation, we next introduce an approximation approach, which leverages low-rank structures and avoids assembling large matrices.

A low-rank approximation $A B \approx C \in \R^{M^d \times M^d}$ tries to capture most of the action of the matrix $C$ by a product of two low-rank matrices $A \in \R^{M^d \times l}$ and $B \in \R^{l \times M^d}$. Several efficient algorithms, e.g., the fast multipole method and $\mathcal{H}$-matrices, depend on low-rank approximations \cite{LH03, GR87}. We approximate eigenvalues and eigenvectors of the correlation matrix $C$ by using low-rank approximations as follows.

Given a low rank approximation $Q (Q' C)$ of the symmetric covariance matrix $C$, where the matrix $Q$ is of size $\R^{M^d \times l}$ with $l \geq D$ orthonormal columns, the eigenvalue problem of  $C$ can be approximated efficiently by applying $QR$ or $SVD$ algorithm to the much smaller matrix $Q' C Q$. In the DgPC setting \eqref{eq:cov_pce}, this amounts to computing
\begin{align} \label{eq:cov_QQ'}
 \R^{l \times l} \ni Q' C Q = \sum_{\balp >0} (Q'u_{\balp}) \, (Q' u_{\balp})' . 
\end{align}
The explicit assembly of the covariance matrix $C$ is avoided by computing only the matrix-vector product $Q' u_{\alpha} \in \R^{l \times 1}$. An approximate eigenvalue decomposition $C \approx U \Lambda U' $ is deduced from the eigenvalue decomposition  of the smaller matrix $Q' C Q = V \Lambda V'$ by setting $U = Q V$.  

The crucial step of the computation is the construction of a low-rank matrix $Q$ with $l \ll M^d$ orthonormal columns that accurately describes $C$.  We tried an approach based on the discrete unitary Fourier transform to map the coefficients $u_{\alpha}$ to the frequency space and retain only the lowest frequencies. Although this approach allowed us to obtain reasonable compressions of the covariance matrix and enabled faster computations, the following approach consistently yielded much better results.

Following \cite{MRM11,HMT11}, we construct the matrix  $Q$ using random projections. Algorithm 4.1 in \cite{HMT11} draws an $M^d \times l$ Gaussian random matrix $O$ and forms the matrix $Y = C O \in \R^{M^d \times l}$. The matrix $Q$ with $l$ orthonormal columns is then obtained by the $QR$ factorization $Y= QR$. Note again that we do not assemble the matrix $C$. Rather, the matrix-matrix product $C O$ is computed as in equation $\eqref{eq:cov_QQ'}$. Since we require the largest $D$ eigenvalues, the target low-rank becomes $D$. As indicated in \cite{HMT11}, we use an oversampling parameter $p$ by setting $l=D+p$. Typical values of $p$ are $5$ or $10$.  Since eigenvalues decay rapidly in our applications, we found $p=10$ to be accurate. With overwhelming probability,  the spectral error $||C - Q Q' C||_2$ is bounded within a small polynomial perturbation of the theoretical minimum, the $(l+1)^{th}$ singular value of $C$; for relevant theoretical details, see \cite[Section 10]{HMT11}.

In practice, we found this randomized approach to be highly accurate in our computational simulations. Moreover, since assembly of large covariance matrices are avoided, running times and memory requirements for the KL expansion in $\R^2$ are reduced drastically compared to the previously described methods; see section \ref{sec:SNS}.

\subsection{Additional non-forcing random inputs} \label{sec:random_viscosity}

In this subsection, we consider the case in which the differential operator $\mathcal{L}$ in \eqref{eq:gen_spde} contains additional random input parameters, i.e., $\mathcal{L}= \mathcal{L}(u(x,t,\omega),\omega)$.  The random inputs will be denoted by the process $Z(x,\omega)$ of a dimension $D_{Z}$. A typical case is that of a random viscosity, e.g., depending on a set of uniformly distributed random variables. We assume that the process $Z$ is independent of time and that the corresponding orthogonal polynomials are available; for instance in the Askey family \cite{XK02}.

We first observe that the solution $u(x,t,\omega)$ is now a functional of Brownian motion $W$ and the random process $Z$. Therefore, assuming $L^2$ integrability, it can be written as a PCE in terms of the associated orthogonal polynomials of $W$ and $Z$. At the restart $t_j$, there are two options to carry out the KL expansion: (i) apply the KL expansion to only the solution $u_j$ and keep the basis variables for $Z$ in addition to $\bta_j$ in the next PCE; and (ii) compress the combined random variable $(u_j,Z)$ using the KL expansion and denote by $\bta_j$ the combined random modes which represent both $u_j$ and $Z$. The first approach will yield PCEs which provide functional representations of the solution in terms of $W$ and $Z$ for each time $t_j$. In the second approach, the functional dependence of the solution in terms of $Z$ is lost in the first KL expansion step. However, the moments of the solution can still be computed through the combined random KL modes. In many UQ settings, rather than a functional dependence, it is statistical information of the underlying solution, e.g., moments of the invariant measure in the long-time, that we are after. Moreover, the second approach can be seen as a dimensionality reduction technique, which compresses $u_j$ and the process $Z$ together, thereby further reducing the number of terms in PCE. When additional random parameters appear in the equation, we found it reasonable to implement the second approach to reduce cost while the first approach may be used as a reference computation to assess accuracy. 

\begin{remark} \rm 
It is useful to note that by combining the random fields $u_j$ and $Z$, the algorithm automatically chooses the important part of the random process $Z$ that influences the solution while keeping the moments of the solution accurate; see section \ref{sec:Numerical}. 
\end{remark}

\subsection{Moments and orthogonal polynomials by a sampling approach}

After obtaining the random modes $\bta_j$, $j \in \mathbb{N} \cup \{ 0 \}$, we need to construct the following  orthonormal basis: 
\begin{align} \label{eq:orth_basis}
\{ T_{\balp}(\bxi_j, \bta_j) \, : |\balp| \leq N, \balp \in \mathcal{J}^r_{K+D,N} \}.
\end{align}
Notice that since $\bxi_j$ is Gaussian and identically distributed for each $j$, the corresponding orthonormal polynomials of $\bxi_j$ are known to be the Hermite polynomials for each $j$. However, the probability distribution of $\bta_j$ is arbitrary and changes at each restart. Therefore, the computation of orthonormal polynomials is computationally intensive and can be performed using the Gram--Schmidt method as follows. 

We note that the set \eqref{eq:orth_basis} can be computed based on the knowledge of moments of the variables $\bxi_j$ and $\bta_j$. Following \cite{GM10,AGPRH12}, we assemble the Gram matrix $H^j$ with the entries
 \begin{align} \label{eq:gram_matrix}
 H^j_{kl} = \E[  (\bxi_j, \bta_j)^{\balp_k+\balp_l}  ], \quad \balp_k,\balp_l \in  \mathcal{J}^r_{K+D,N}. 
 \end{align}
 The matrix $H^j$ is a $|\mathcal{J}^r_{K+D,N}|$-dimensional, square and symmetric matrix. For theoretical reasons, we assume that the moments up to $2N$ exist and the measure of $(\bxi_j,\bta_j)$ is non-degenerate. Then, the Cholesky factorization is employed to $H^j$ and the polynomials $T_{\balp_l}$ is found by inverting the resulting upper triangular matrix as
\begin{align} \label{eq:Talp}
T_{\balp_l}(\bxi_j, \bta_j) = \sum_{\balp_k \leq \balp_l}  a_{kl}  \, (\bxi_j, \bta_j)^{\balp_k}, 
\end{align}
where $a_{kl}$ are real coefficients.

\begin{remark} \rm 
The KL expansion yields uncorrelated random variables $\bta_j$. If the underlying process is Gaussian, it is known that these variables are also independent. However, in general, marginals of $\bta_j$ are dependent variables. Multi-index operations can still be used to construct the polynomial set \eqref{eq:orth_basis} with respect to the joint distribution, although the estimation of multivariate moments of $\bta_j$ becomes necessary because of such a dependency. In this case, it is known that orthogonal polynomials are not unique and depend on the ordering imposed on the multi-index set; see \cite{OB16,AGPRH12,EMSU12}. In all computations, we use the graded lexicographic ordering for multi-indices. 
\end{remark}

\begin{remark} \rm 
The completeness of the orthogonal polynomials $T_{\balp}(\bxi_j, \bta_j)$ is closely related to the moment problem of the random variables $\bxi_j$ and $\bta_j$. In particular, if the moment problem is uniquely solvable, i.e., the measure is determinate, then the orthogonal polynomials are dense in $L^2$  ~\cite{BC81,P82, gautschi, EMSU12, OB16}. Some basic conditions that guarantee determinacy of the measure of a continuous random variable on a finite dimensional space are compact support and exponential integrability.  Gaussian measures are determinate and the Hermite PCE converges by the Cameron-Martin theorem. However, in general, whether the distribution of $\bta_j$ is determinate or not is unknown. This problem is addressed in our previous work in the case of finite dimensional SDE systems; see \cite[section 3]{OB16}. Theoretical results are applied in the setting where the solutions are approximated by compactly supported distributions under appropriate assumptions. In the following, we assume that the measures associated to $\bta_j$ are determinate so that convergence is guaranteed, which is consistent with our numerical simulations; see section \ref{sec:Numerical}.
\end{remark}

Based on the above discussion, the orthonormal basis \eqref{eq:orth_basis} requires the computation of the moments \eqref{eq:gram_matrix}. The exact moments of the Gaussian variables $\bxi_j$ are computed by analytical formulas and then stored during the offline stage. However, the distribution of $\bta_j$ is varying with $j$. Therefore, the computation of moments should be carried out based on information provided by the PCE. 

Several methods are available to compute moments of probability distributions in the PC context such as, e.g., quadrature methods, Monte Carlo sampling, or a pure PCE approach. This procedure is notoriously expensive and ill-posed \cite{gautschi}. The pure PCE approach computes the moments of $\bta_j$ by repeatedly multiplying its PCE and taking expectation; see \cite{DNPKGL04,OB16}. This approach is discussed in detail in our earlier manuscript \cite{OB16} and works reasonably well for a low dimensional SDE systems. However, it becomes prohibitively expensive if the dimension of the random variables in the PC basis is even moderate; see the computational complexity section in \cite{OB16}. Therefore, in this work, we consider an alternative approach using Monte Carlo sampling, which drastically reduces the computational cost for computing moments compared to the pure PCE approach.    

We assume that independent samples of the initial condition (therefore, the samples of $\bta_0$) are provided so that the algorithm can be initialized. To construct the set \eqref{eq:orth_basis}, based on \eqref{eq:gram_matrix},  we need to compute the first $2N$ moments of the joint random variable $(\bxi_j,\bta_j)$. Moreover, since the triple products will be required for the evolution of the PCE coefficients, the first $3N$ moments need to be computed as well; see section \ref{sec:galerkin_proj}. Using the same ordering of the multi-index set $\mathcal{J}^r_{K+D,3N}$ , we require the computation of the following moments: 
\begin{align} \label{eq:mom_comb}
\E[ (\bxi_j, \bta_j)^{\balp} ] =  \E [\bxi_j^{(\alpha_1,\ldots,\alpha_K)}] \,\E [\bta_j^{(\alpha_{K+1},\ldots,\alpha_{K+D})}]   ,\quad \balp \in \mathcal{J}^r_{K+D,3N},
\end{align}
where we used the independence of $\bxi_j$ and $\bta_j$.

Let $\bta_j (\omega_i) :=(\eta_{j,1}(\omega_i),\ldots,\eta_{j,D}(\omega_i)) $ denote independent samples of the random modes for $\omega_i \in \Omega$, where $i=1,\ldots,S \in \mathbb{N}$. Then, provided the samples $\bta_j(\omega_i)$ are given, the moments of $\bta_{j}$ can be approximated by
\[
 \E [\bta_j^{(\alpha_{K+1},\ldots,\alpha_{K+D})}]  \approx \frac{1}{S} \sum_{i=1}^S (\bta_j(\omega_i))^{(\alpha_{K+1},\ldots,\alpha_{K+D})}  =\frac{1}{S} \sum_{i=1}^S  \prod_{l=1}^D (\eta_{j,l}(\omega_i) )^{\alpha_{K+l}},
\]
where we used the usual multi-index notation for powers. Therefore, multivariate moments \eqref{eq:mom_comb} are computed by a combination of the analytical formulas  for $\bxi_j$ and a sampling approximation for $\bta_j$. Note that in applications, we use small values of $N$ with a sufficiently large number of samples $S$ to guarantee accuracy. 

Although we discussed computing moments based on samples, we have not explained how the samples of $\bta_j$ are acquired except for $\bta_0$. 
The distribution of $\bta_j$, $j \geq 1$, is evolving in time. However, owing to the PCE \eqref{eq:etaij} of the each component $\eta_{j,l}$ , we can write
\[
\bta_j = \sum_{\balp} \bta_{j,\balp} \, T_{\balp}(\bxi_{j-1}, \bta_{j-1}),  \quad j \in \mathbb{N}.
\]
This representation gives a natural way to sample from the distribution of $\bta_j$ by the recursion
\begin{align} \label{eq:sample_diagonal}
\bta_j(\omega_i) = \sum_{\balp} \bta_{j,\balp} \, T_{\balp}(\bxi_{j-1}(\omega_i), \bta_{j-1}(\omega_i)) , \quad i=1,\ldots,S,
\end{align}
assuming that we obtained samples of $\bta_{j-1}(\omega_i)$ at the previous restart $t_{j-1}$. Independent samples $\bxi_{j-1}(\omega_i)$ are obtained through sampling a multivariate Gaussian distribution. Therefore, on the subinterval $[t_{j-1},t_j]$, PCE acts like a transport map which maps previously obtained samples of $\bta_{j-1}$ and new samples of $\bxi_{j-1}$ to the samples of the new random modes $\bta_j$. 

\begin{remark} \rm Note that $\bta_j$ is a function of the variables $\bxi_{j-1}$ and $\bta_{j-1}$. Therefore, the number of samples of $\bta_j$ should be ideally $S^2$ provided the same of number of samples $S$ is used for each $\bxi_{j-1}$ and $\bta_{j-1}$. However, in practice, this is not feasible in our method as the number of samples grows in time. Instead, the equation \eqref{eq:sample_diagonal} keeps only the diagonal terms in the sample space. In the numerical simulations, we use large values of $S$ so that the loss of accuracy incurred from discarding some samples would be minimal. Moreover, it is important to notice that \eqref{eq:sample_diagonal} entails samples from the joint distribution of $\bta_j$ so that Monte Carlo method is used to approximate the expectations while preserving the dependence structure of marginals.
\end{remark}

\begin{remark}\label{rmk:samples} \rm 
The method blends Monte Carlo sampling into a PC approach to exploit the virtues of the both methods, namely, rapid computation of expectations and spectral accuracy provided by the MC and PC methods, respectively. Although the method is utilizing samples for the computation of moments, samples are not used in the evolution stage. The algorithm essentially propagates moments of the measures between successive times, where moments are computed using a sampling technique. To test the robustness of the method with respect to sampling, imagine that the algorithm starts with an infinite supply of independent samples of $\bta_{0}$. We discard the first sample after using it to construct the corresponding orthogonal polynomials at the end of the first time interval and propagate the remaining samples with the PCE map to construct (an infinite supply of) samples of $\bta_1$. The algorithm iteratively estimates the distribution, hence the moments, of each $\bta_j$ while samples are discarded at each restart. We tested this idea by starting with a set of $n$ independent samples of $\bta_0$ and propagated them by PCE for a maximum of $n$ restarts. We found that the accuracy of the calculations was not affected by such a re-sampling tool. This comparison showed the robustness of the algorithm under changes of samples. Since in practice, such re-sampling increases the computational costs compared to \eqref{eq:sample_diagonal}, it is not considered in the numerical experiments presented in the next sections. We also emphasize that the sampling approach readily returns samples of the approximated solution at the endpoint $T$ through its KL expansion without a further sampling procedure. These samples can also be useful in uncertainty quantification to  estimate further statistical properties of the solution such as probabilities on prescribed sets or probability density functions. 
\end{remark}

\subsection{Galerkin projection and local initial conditions} \label{sec:galerkin_proj}

Once an orthonormal basis is obtained, the algorithm performs a Galerkin projection onto the space spanned by the basis, and this requires the computation of triple products.  Using the representation \eqref{eq:Talp}, the required triple products can be written as  
\begin{align} \label{eq:tripprod_formula}
\E[ T_{\balp_{k}} T_{\balp_{l}} T_{\balp_m}] = \sum_{ \balp_{k'} \leq \balp_k}\sum_{ \balp_{l'} \leq \balp_l} \sum_{ \balp_{m'} \leq \balp_m}  a_{k'k} \, a_{l'l} \, a_{m'm} \, \E[(\bxi_j,\bta_j)^{\balp_{k'} +\balp_{l'}+ \balp_{m'} }],
\end{align}
where all multi-indices belong to the set $\mathcal{J}^r_{K+D,N}$. Thus this formula can be computed by the knowledge of moments of order up to $3N$. 

Depending on the choice of the sparse index $r$, the multi-index $\balp_{k'} +\balp_{l'}+ \balp_{m'} \in \mathcal{J}_{K+D,3N}$ might not be an element of $\mathcal{J}^r_{K+D,3N}$. Therefore, once we fix the multi-index set $\mathcal{J}^r_{K+D,N}$ in the offline stage, we also compute and store the indices that are elements of $\mathcal{J}^r_{K+D,3N}$. 

Finally, we perform Galerkin projection of the SPDE \eqref{eq:gen_spde} and obtain the following PDE system for the coefficients $u_{j+1,\balp}(x,t)$ of \eqref{eq:PC_exp}, $t \in [t_j,t_{j+1}]$ : 
\begin{align} \label{eq:PDE_sys}
\partial_t (u_{j+1,\balp}) = \E \left[T_{\alpha}(\bxi_j, \bta_j) \, \mathcal{L}\left(\sum_{\boldsymbol{\beta}} u_{j+1,\boldsymbol{\beta}} \, T_{\boldsymbol{\beta}}(\bxi_j, \bta_j) \right) \right] + 
\sigma(x) \, \E [T_{\balp}(\bxi_j, \bta_j) \dot{W}(t)].
\end{align}
The first expectation in the above line is computed with the aid of the triple products \eqref{eq:tripprod_formula} and the second using the representation \eqref{eq:conv_W}.

Note that the initial conditions $u_{j+1,\balp}(x,t_j)$ can be obtained by noticing that the representation \eqref{eq:KL_PC} of $u_j$ is nothing but a sum involving linear polynomials in $\eta_{j,l}$. It can therefore be rewritten in the basis $T_{\alpha}(\bxi_j, \bta_j)$ with the help of Galerkin projections. Hence, the only coefficients that survive in \eqref{eq:PC_exp} at $t=t_j$ are the mean and the ones which correspond to the first degree polynomials in $\bta_j$. Then, the PDE system \eqref{eq:PDE_sys} can be solved in time using a time-integration method combined with the aforementioned spectral method \eqref{eq:cov_spectral}. If the initial condition $u_0$ of \eqref{eq:gen_spde} is deterministic, we employ the Hermite PCE on the first interval $[0,t_1]$, which does not necessitate the computation of the KL decomposition. 
 
\subsection{Adaptive restart scheme}

So far, the method uses a predetermined restart time $\Delta t$. For long-term simulations, an adaptive restart scheme that sets the restarts online depending on the properties of the solution can reduce the computational cost. 

We propose to adapt the restart time based on the following two criteria: (i) preserve the accuracy of the representation \eqref{eq:conv_W} of the random forcing;  and (ii) mitigate the effect the nonlinearities in the accuracy of the polynomial expansions. For a prescribed number of dimensions to describe the random forcing, the algorithm can not take too large steps to preserve accuracy. Also, nonlinearities force the solution to be less accurately described by low-degree polynomials in the initial condition as time increases. In both cases, we wish $\Delta t$ to be as large as possible for a given accuracy in mind.

To this end, let $\Delta  t_j$ denote the adaptive time-step starting from time $t_j$. To ensure an accurate representation of the forcing term, we set a maximum value $\Delta t_{\max}$ for $\Delta t_j$ for all $j$, i.e. $\Delta t_j \leq \Delta t_{\max}$. In practice, $\Delta t_{\max}$ is based on the error analysis of random forcing by a finite dimensional approximation; see for instance \cite{HLRZ06}. To address nonlinearities, we consider the following ratio for the PC coefficients $u_{j,\balp}$: 
\begin{align} \label{eq:rho}
 \rho(t):= \frac{|| \sum_{|\balp|>1} u_{j,\balp}^2(\cdot,t) ||_{L^1}}{|| \sum_{|\balp|>0} u_{j,\balp}^2(\cdot,t) ||_{L^1}}, \quad t \in [t_{j-1}, t_{j-1}+ \Delta t_{j-1}].
\end{align}
The condition measures the norm ratio of the nonlinear terms in the variance to the norm of the variance. In applications, the ratio is computed at each time integration point. Similar other conditions were used in different settings in \cite{GSVK10, HS14}. 

Consider a threshold value $\epsilon \in (0,1)$. We propose the following conditions for adaptive time-steps using $t \in [t_{j-1}, t_{j-1}+ \Delta t_{j-1}]$ :
\begin{enumerate}[i)]
\item $
\mbox{if  } \rho(t) \leq 3 \epsilon  \mbox{  then  set the next time-step } \Delta t_j =  {\rm min}(t_* - t_{j-1},\Delta t_{\rm max}) \mbox{ for } \rho_p(t_*)=2\epsilon $, \label{item:less3eps}
\item $\mbox{if  } \rho(t) >3 \epsilon  \mbox{  then  go back to } t_{j-1} \mbox{ and set } \Delta t_{j-1} =   {\rm min}(t_* - t_{j-1},\Delta t_{\rm max}) \mbox{ for } \rho(t_*)=2\epsilon $, \label{item:bigger3eps}
\end{enumerate}
where $\rho_p$ is a polynomial approximation to $\rho(t)$, which can be found by fitting a $p$-degree polynomial to $\rho(t)$ on the interval  $ [t_{j-1}, t_{j-1}+ \Delta t_{j-1}]$. This approximation is only required for time values $t$ satisfying $\rho(t) < 2\epsilon$. 

The time-steps $\Delta t_j$, $j >0$ are set adaptively. For short time-steps, we do not expect dynamics to change drastically between successive intervals. Therefore, condition \ref{item:less3eps} verifies whether the ratio is smaller than $3 \epsilon$ on the current interval, and then sets the adaptive time-step for the next interval. When $\rho(t) \leq \epsilon$, then the algorithm selects a bigger time-step by finding the root of $\rho_p(t_*)=2\epsilon$. Note that the current evolution on the interval $[t_{j-1}, t_{j-1}+ \Delta t_{j-1}]$ is not prematurely stopped at the end point. Although PCEs converge at any point inside the interval, errors, however, are known to wildly oscillate inside the interval and become spectrally accurate only at the end point; see \cite{BM13, OB16}. (This is because $W(t)$ in \eqref{eq:conv_W} is spectrally accurate at $t=T$ and much less so on $(0,T)$.)  Condition \ref{item:bigger3eps} essentially verifies whether the ratio becomes too large (i.e. $>3 \epsilon$), and when this happens, forces the evolution to restart from the current initial point $t_{j-1}$. This control ensures that the algorithm does not take too large steps.

Our procedure is summarized in Algorithm \ref{alg:DgPC}, where, for simplicity, we only present the version which uses a predetermined number of restarts.

\begin{algorithm}
\caption{Dynamical generalized Polynomial Chaos (DgPC) for SPDEs}
\label{alg:DgPC}
\begin{algorithmic}

\\ Decompose the time domain $[0,T] = [0,t_1] \cup \ldots \cup [t_{n-1},T]$
\\ Initialize the degrees of freedom $K,N,D,S$
\\ Choose the sparse index $r=(r_1,\ldots,r_{K+D})$
\\ Compute the indices used in the triple-product formula \eqref{eq:tripprod_formula}
\\ Compute moments of $\bxi_0$
\For{each time-step $t_j$}
\State apply the KL expansion to $u_j$ and obtain $\bta_j=(\eta_{j,1},\ldots,\eta_{j,D})$
\State compute the moments $\E [(\bxi_j, \bta_j)^{\balp}]$
\State construct orthogonal polynomials $T_{\balp}(\bxi_j, \bta_j)$
\State compute the associated triple products 
\State perform Galerkin projection onto $\mbox{span} \{ T_{\balp} (\bxi_j,\bta_j) \}$
\State set up initial conditions for $u_{j+1}$
\State evolve the PCE coefficients of $u_{j+1}$
\EndFor 
\end{algorithmic}
\end{algorithm}

\section{Numerical simulations} \label{sec:Numerical}

We now present numerical simulations for the Burgers equation in one spatial dimension and a Navier--Stokes system in two spatial dimensions both driven by white noise. We consider these equations for two reasons. First, the statistical behavior of solutions of these equations is of importance in statistical mechanics and turbulence theory; see, e.g., \cite{GS91,MR04, HM06, HM08}.  Second, they serve as challenging test beds for the PCE methodology. 

We illustrate convergence results in terms of the degrees of freedoms $D$, $N$ and the time-step $\Delta t$, and consider both short-time and long-time evolutions. The convergence of the method in terms of $K$ is treated in detail in \cite{OB16}. As a general comment, we do not recommend using large values of $N$ since the computation of orthogonal polynomials is quite ill-posed. In the following, we use polynomials of degree up to $N=3$. The algorithm mitigates the ill-posedness by choosing frequent restarts and small degree; i.e. small $\Delta t$ and $N$.  The settings we consider here closely follow those addressed in the manuscripts \cite{HLRZ06,Luo}.

\subsection{Burgers equation}

We consider the following one dimensional viscous stochastic Burgers equation for $u(x,t,\omega)$: 
\begin{align} \label{eq:1DBurgers}
\begin{cases}
 \partial_t u + \frac{1}{2}  \partial_x u^2 = \nu \, \partial^2_x  u  + \sigma(x) \dot{W}(t), \\
 u(x,0,\cdot) = u_0(x) , \quad u(0,t,\cdot) = u(1,t,\cdot), \quad (x,t) \in [0,1] \times [0,T],
 \end{cases}
\end{align} 
where $W(t)$ is a Brownian motion in time, $\nu>0$ is the viscosity (which will be either deterministic or random), the initial condition $u_0$ is deterministic, and the solution itself is periodic in the spatial domain $[0,1]$. 

Following \cite{Luo, HLRZ06, OB16}, we choose cosine functions as an orthonormal basis $m_{j,i}(t), i \in \mathbb{N}$, for $L^2[t_{j},t_{j+1}] $. Employing the equation \eqref{eq:conv_W} for each subinterval $[t_j,t_{j+1}]$ and using Galerkin projection, we derive the governing equations for the PC coefficients $u_{j+1,\balp}(x,t)$ of $u_{j+1}$: 
\[
\partial_t (u_{j+1,\balp}) + \frac{1}{2} \partial_x \,  (u^2_{j+1})_{\balp} = \nu \, \partial^2_x \, u_{j+1,\balp} + \sigma(x) \sum_{i=1}^K m_{j,i}(t) \, \E [   \xi_{j,i} \, T_{\balp}(\bxi_j,\bta_j)  ],
\] 
where $\bxi_j = (\xi_{j,1}, \xi_{j,2},\ldots)$. Since the initial condition is deterministic, we employ Hermite PCE in the subinterval $[t_0,t_1]$. The PC coefficients  $(u^2_{j+1})_{\balp}$ of the nonlinearity $u^2_{j+1}$ are computed by multiplying the corresponding PCEs with the help of pre-computed triple products. 

Since the coefficients $u_{j+1,\balp}(x,t)$ are periodic in the physical space, we utilize a truncated Fourier series:
\[
u_{j+1, \balp}(x,t) \approx \sum_{k=-M/2+1}^{M/2} \hat{u}_{j+1,\balp}(k,t)\,  e^{2\pi \, i \, k \, x}, \quad x \in [0,1], 
\]
with the even number of frequencies $M$ to be chosen. Further, using the equidistant partition for the spatial domain $[0,1]$ and Fast Fourier Transform (FFT), we can compute the Fourier coefficients  
 $\hat{u}_{j+1,\balp}$ with a reasonable computational cost. This procedure gives rise to an ODE system which is then integrated in time using a second order predictor-corrector method combined with an exponential integrator for the stiff term~\cite{L07}. 

In the implementation of the Burgers equation, the algorithm assembles the covariance matrix at each restart as discussed in section \ref{sec:KLexpansion}. A large memory is not required for such a one dimensional spatial problem. We found that using the random projection technique described in section \ref{sec:KLexpansion} did not result in significantly shorter total computational times because the computation of the eigenvalue problem is already efficient in this case by means of Krylov subspace methods.

\begin{exmp}\label{ex:Burgers1} \rm

For this numerical simulation,  we choose the spatial part of the random forcing as $\sigma(x)= \frac{1}{2}\cos(2 \pi x)$, the initial condition  as $u_0(x) = \frac{1}{2} \sin(2 \pi x)$ and set the viscosity $\nu=0.01$. Under these sets of parameters, it has been proved that there exists a unique invariant measure, which is the long-time limit of the dynamics \cite[Theorem 2]{S91}. Thus, the main aim of the following simulations is to demonstrate the efficiency of the algorithm in the long-time setting.
 
For the parameters of DgPC, we take $K=2$, $N=2$, and vary the number of the KL modes $D$. Final time is $T=3$ and the interval divided into $30$ pieces by taking $\Delta t =0.1$.   The spatial mesh size $M$ is set to be $2^7$ and the number of samples $S$ to compute moments is taken as $10^5$. Finally, sparse indices $r$ are chosen as follows: 

\begin{enumerate}[i)]
\item if $D=3$, then $r=(2,2,2,1,1)$, and if $|\balp| =2$, we set $r=(2,2,2,1,\cdot)$ resulting in $15$ terms in the expansion.
\item if $D=4$, then $r=(2,2,2,2,1,1)$, and if $|\balp| =2$, we set $r=(2,2,2,2,1,\cdot)$ resulting in $21$ terms in the expansion.
\item if $D=5$, then $r=(2,2,2,2,2,1,1)$, and if $|\balp| =2$, we set $r=(2,2,2,2,2,1,\cdot)$ resulting in $28$ terms in the expansion.
\end{enumerate}

Note that the first $K$ indices in $r$ correspond to degrees of polynomials in $\bxi$ and the remaining to $\bta$. Choosing $r=(2,2,2,2,1,1)$ means using only first degree polynomials is $\eta_3$ and $\eta_4$. Also, setting $r=(2,2,2,2,1,\cdot)$ for $|\balp|=2$ eliminates the cross terms involving first degree polynomials of $\eta_4$ in the second order terms. Note that using a sparse index not only reduces the number of terms in PCE but also alleviates the computation of moments. 

To compare our algorithm, we use a second order weak Runge-Kutta scheme since the exact solution is not available. To make a fair comparison, both algorithms (Monte Carlo \& DgPC) use the same time-step $dt=0.001$ and the same mesh size $M=2^7$. The number of samples used in MC algorithms are $N_{\mbox{samp}}=10^4,5\times 10^4$ and $10^5$ and the corresponding algorithms will be denoted by MC1, MC2 and MC3, respectively. The exact solution is taken as the result of MC3 and the relative $L^2$ error $|| \E[u_{dgpc}] - \E [u_{mc}] ||_{2}/|| \E [u_{mc}]||_2 $ of the mean is computed. Errors for higher order centered moments are computed similarly. 

\captionsetup[subfigure]{labelformat=empty}
\begin{figure}[!htb]
\centering
\subfloat[]{
\includegraphics[scale=0.275]{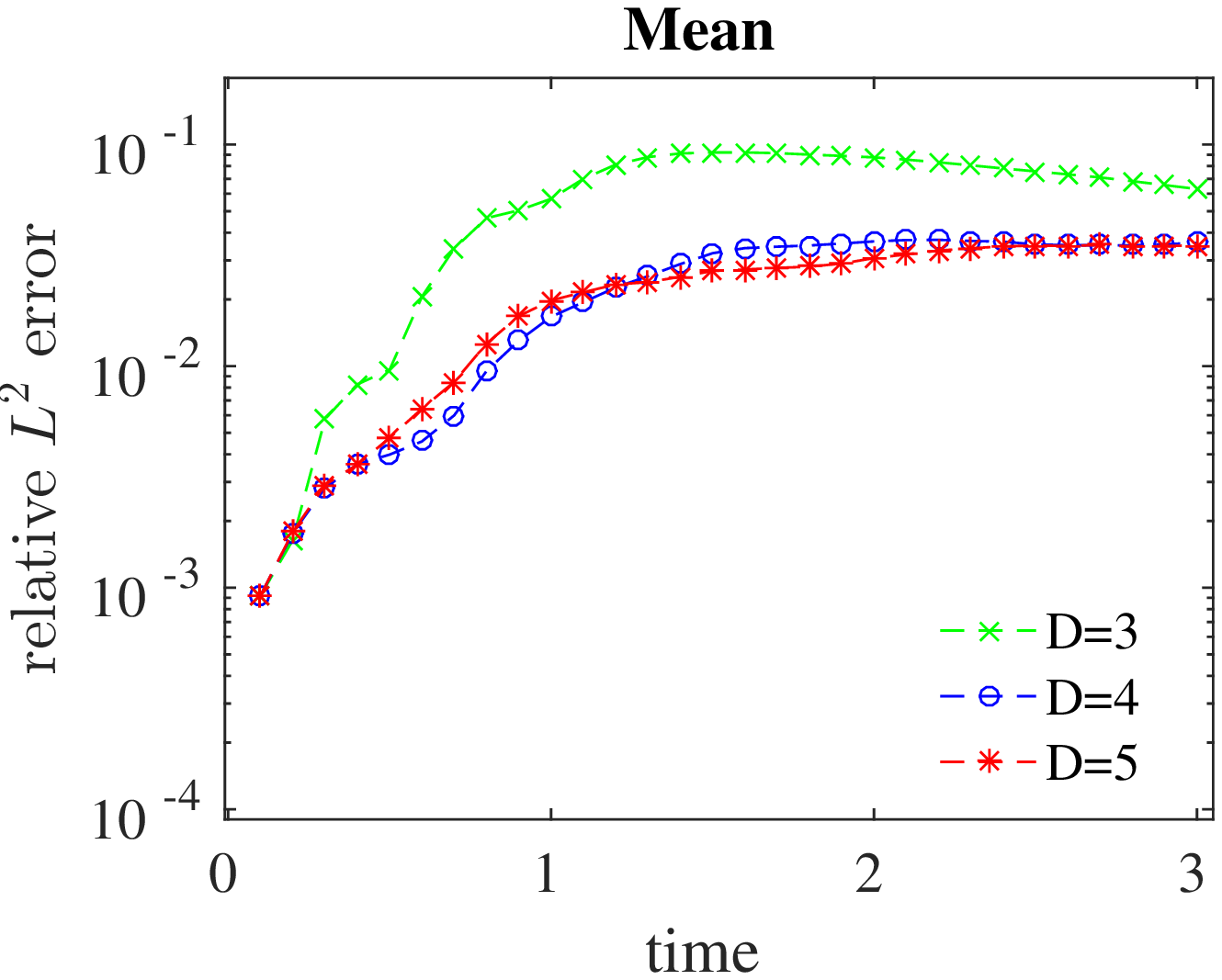}} 
\subfloat[]{
\includegraphics[scale=0.275]{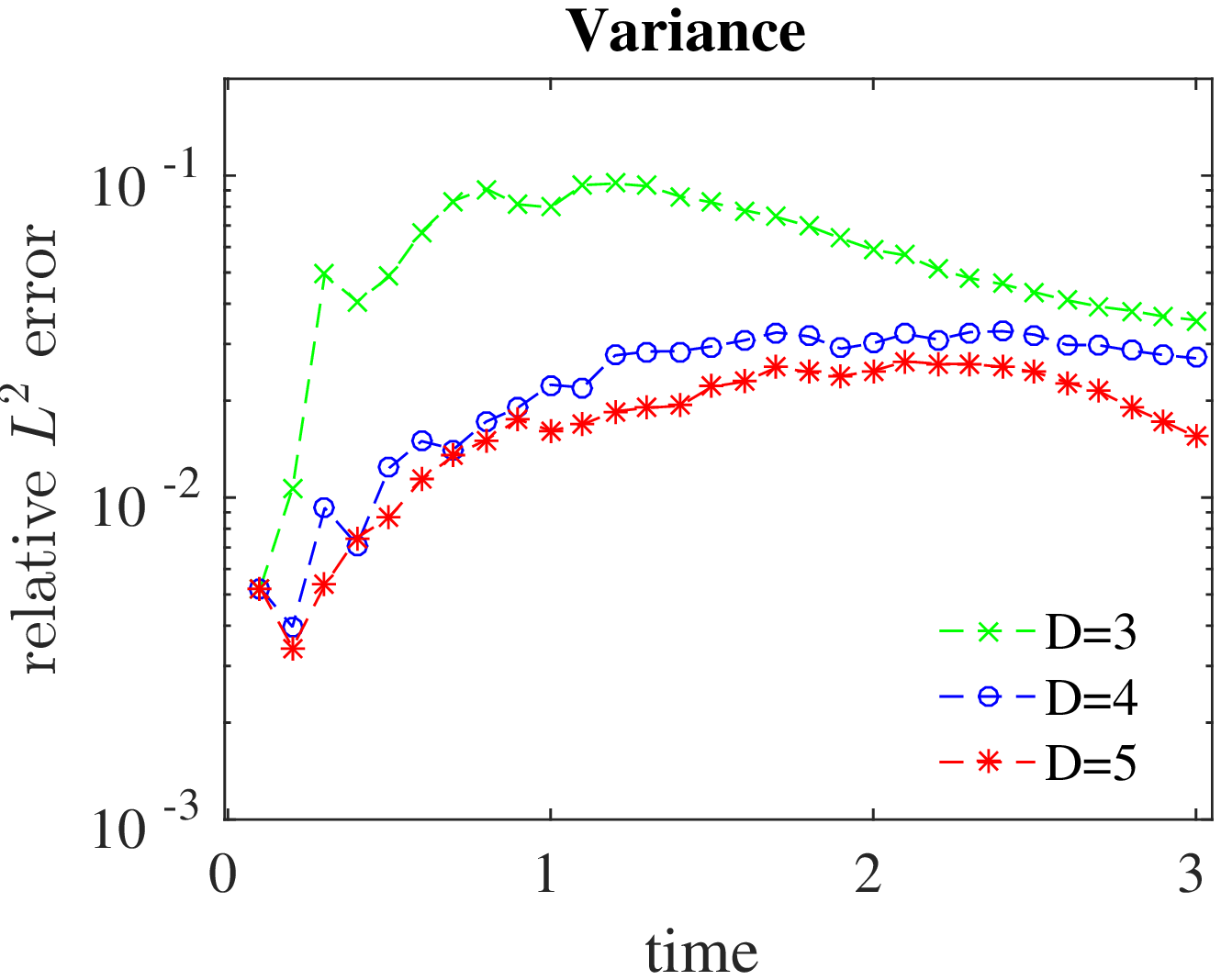}}  
\subfloat[]{
\includegraphics[scale=0.275]{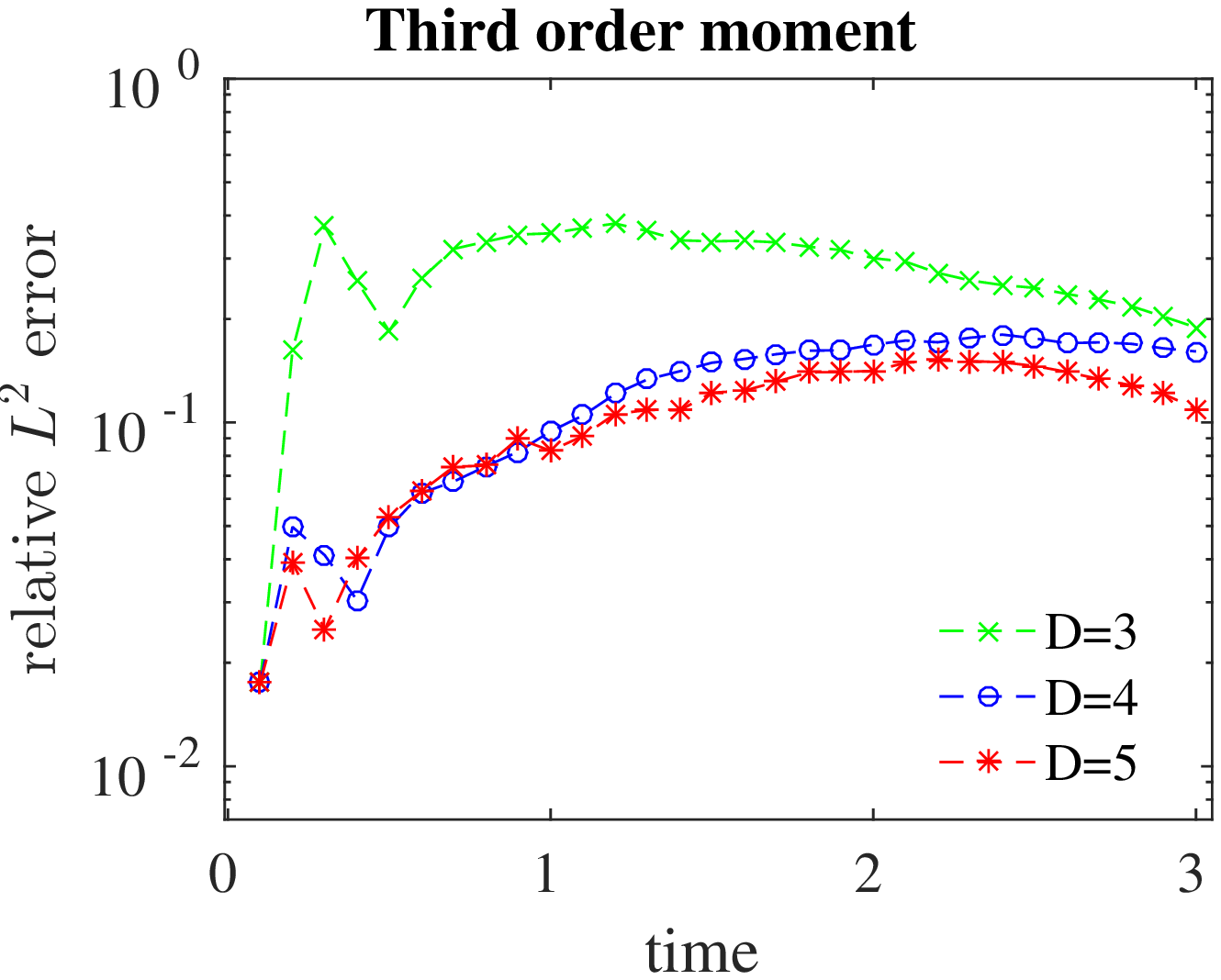}} 
\subfloat[]{
\includegraphics[scale=0.275]{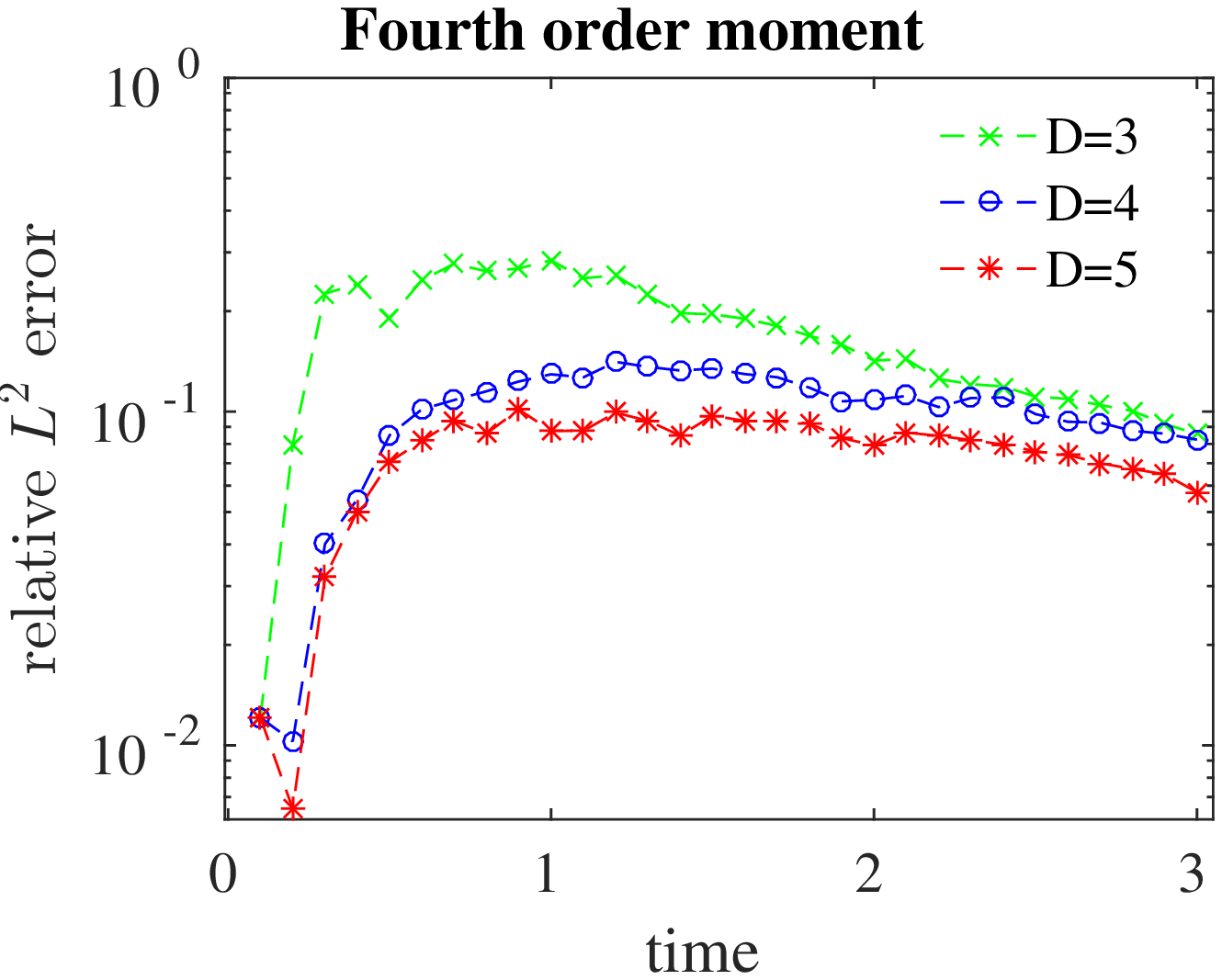}}  
\caption{Relative errors for the centered moments obtained by DgPC with $T=3$ and $\Delta t=0.1$. Exact solution is computed by MC3.}
\label{fig:Ex1_1}
\end{figure} 

From Figure \ref{fig:Ex1_1}, we observe that all errors grow with time in the initial stages and in particular, the degree of freedom $D=3$ is the least accurate, which is expected as the dynamics change rapidly during initial stages. Increasing the number of KL modes entails more accurate expansion up to some order. It can be observed that all the error levels stabilize for moderate times while $D=5$ is the most accurate. This phenomenon is explained by the convergence to a stationary measure such that statistics do not change considerably after some time.   

We now increase the final time to $T=6$. DgPC algorithms use the following parameters: $K=4$, $N=2$, $D=3,4,5$, $S=10^5$  and $\Delta t =0.1$. The corresponding total number of terms for each subinterval becomes $18$, $24$ and $31$. The mesh size is taken as $M=2^8$, which offers better spatial resolution. 
Table \ref{table:burgdet} summarizes the relative $L^2$ errors of the DgPC algorithms with different degrees of freedom and the MC methods. All errors are computed by taking MC3 as the exact solution. The time ratio column is computed as the total time required by the each algorithm divided by the elapsed time of MC3 with $N_{\mbox{samp}}=10^5$ and $dt=0.0005$. The parameters for MC1 and MC2 algorithms remain the same as above; the algorithms are executed a few times and the resulting errors are averaged. We also include the elapsed times for the offline computation in DgPC algorithms. In practice, the required data for the offline step can be computed once and stored for further executions of the algorithm to speed-up the running time. 

\begin{table}[H] 
\centering
\begin{tabular}{r|c|c|c| c |c}
& Mean  & Variance & $3${rd} order & $4${th} order& Time ratio \\ \hline
DgPC: $D=3$ & 2.21E-2 & 1.18E-2 & 5.20E-2& 4.04E-2 & 0.003 \\ 
DgPC: $D=4$ & 1.86E-2& 5.4E-3 & 3.45E-2 & 2.41E-2 & 0.007 \\ 
DgPC: $D=5$ & 1.67E-2 & 4.0E-3& 1.43E-2 & 1.09E-2 & 0.02 \\ 
MC1 & 2.29E-2 & 1.16E-2 & 4.53E-2 & 2.27E-2 &0.05 \\
MC2 & 1.17E-2 & 4.0E-3 & 1.82E-2 & 9.4E-3 &0.25 \\
\hline
\end{tabular}
\caption{Relative errors for the centered moments by DgPC and MC methods at T=6. Each time ratio is computed by comparing to MC3.}
\label{table:burgdet} 
\end{table}

Table \ref{table:burgdet} demonstrates the idea of using low degree polynomials and small number of terms in PC expansion combined with frequent restarts seems to pay off. DgPC with $31$ terms in the expansion (i.e. $D=5$) attains comparable accuracy as MC2 (i.e. $N_{\mbox{samp}}=5 \times 10^4$ and $dt=0.001$) with a computational time which is  only eight percent of that Monte Carlo algorithm. Also, we observe that all errors recovered to a level of $O(10^{-2})$, which is an acceptable accuracy for long-time simulations.

Fixing the degree of freedom as $D=5$, we now employ the adaptive time stepping \eqref{eq:rho} approach. To probe the sensitivity of the algorithm on the threshold parameter $\epsilon$, we choose $\epsilon = 0.005, 0.01, 0.02$  and the initial time-step $\Delta t_0 =0.1$. Also, we set $\Delta t_{max}=0.4$ to get $O(10^{-2})$ accuracy  for the truncation \eqref{eq:conv_W}  of the forcing term  using $K=4$; see \cite[Theorem 5.1]{HLRZ06}. We utilize quadratic polynomials as our ansatz to approximate the ratio \eqref{eq:rho}; see condition \ref{item:less3eps} below \eqref{eq:rho}.

\captionsetup[subfigure]{labelformat=empty}
\begin{figure}[!htb]
\centering
\subfloat[]{
\includegraphics[scale=0.275]{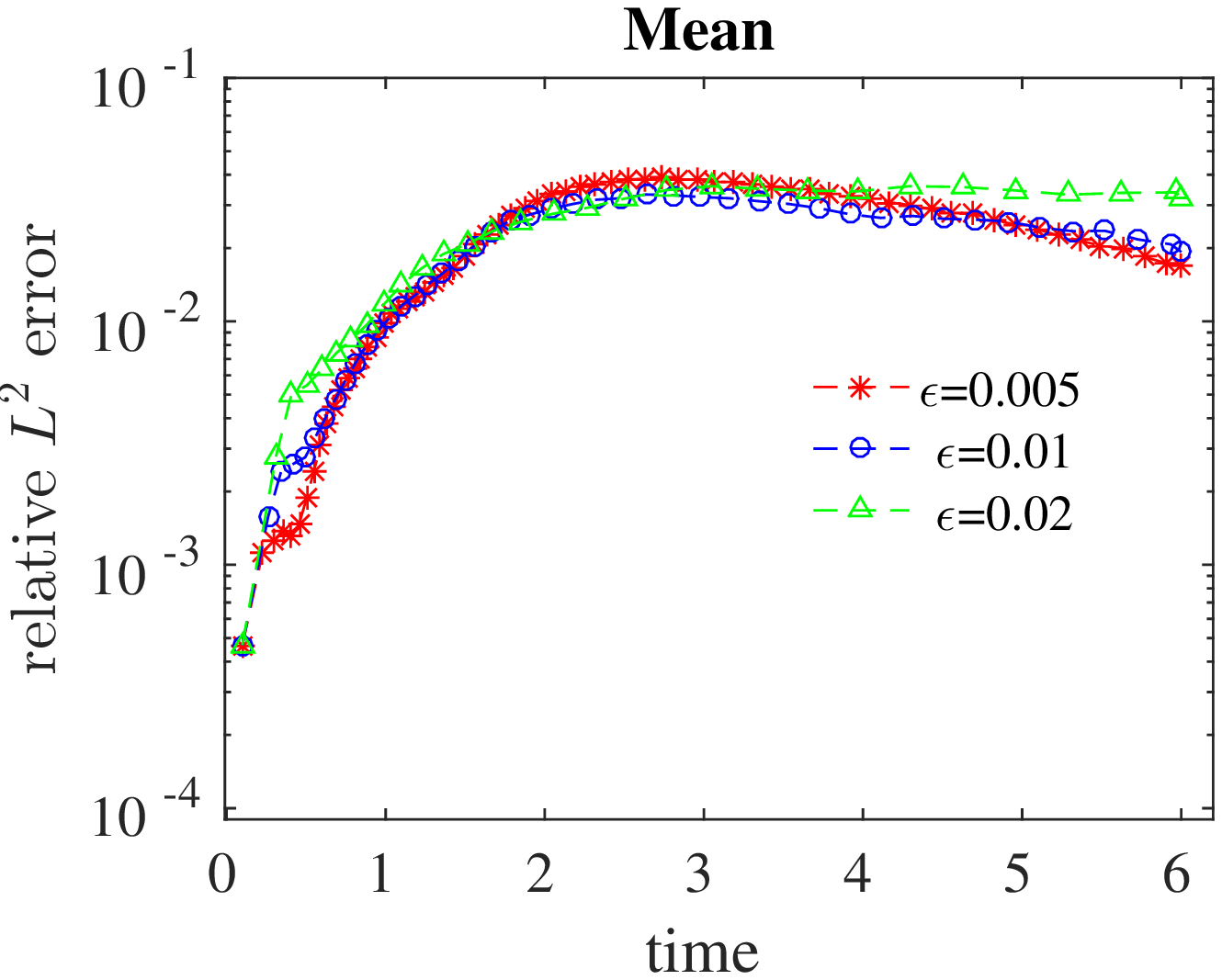}} 
\subfloat[]{
\includegraphics[scale=0.275]{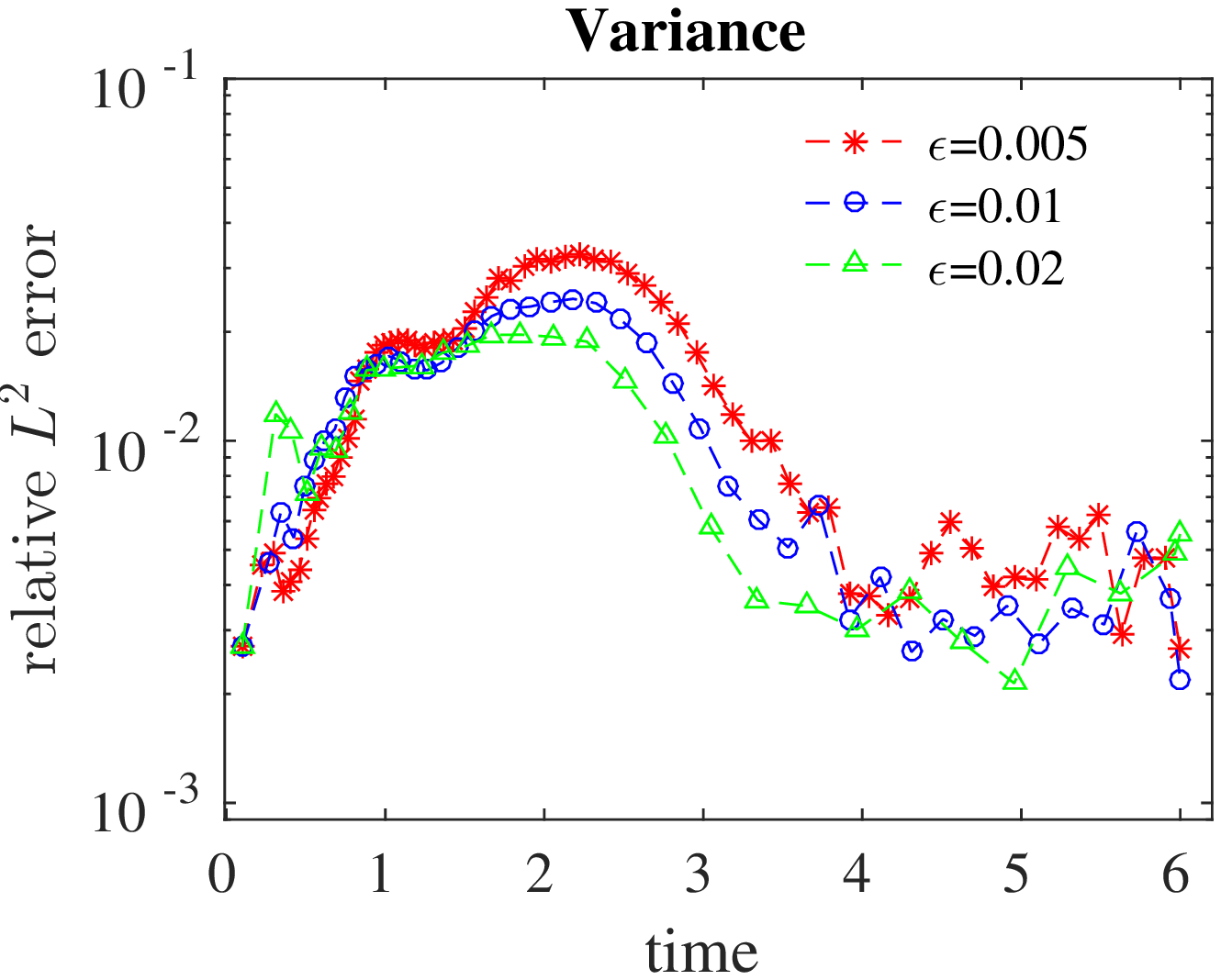}}  
\subfloat[]{
\includegraphics[scale=0.275]{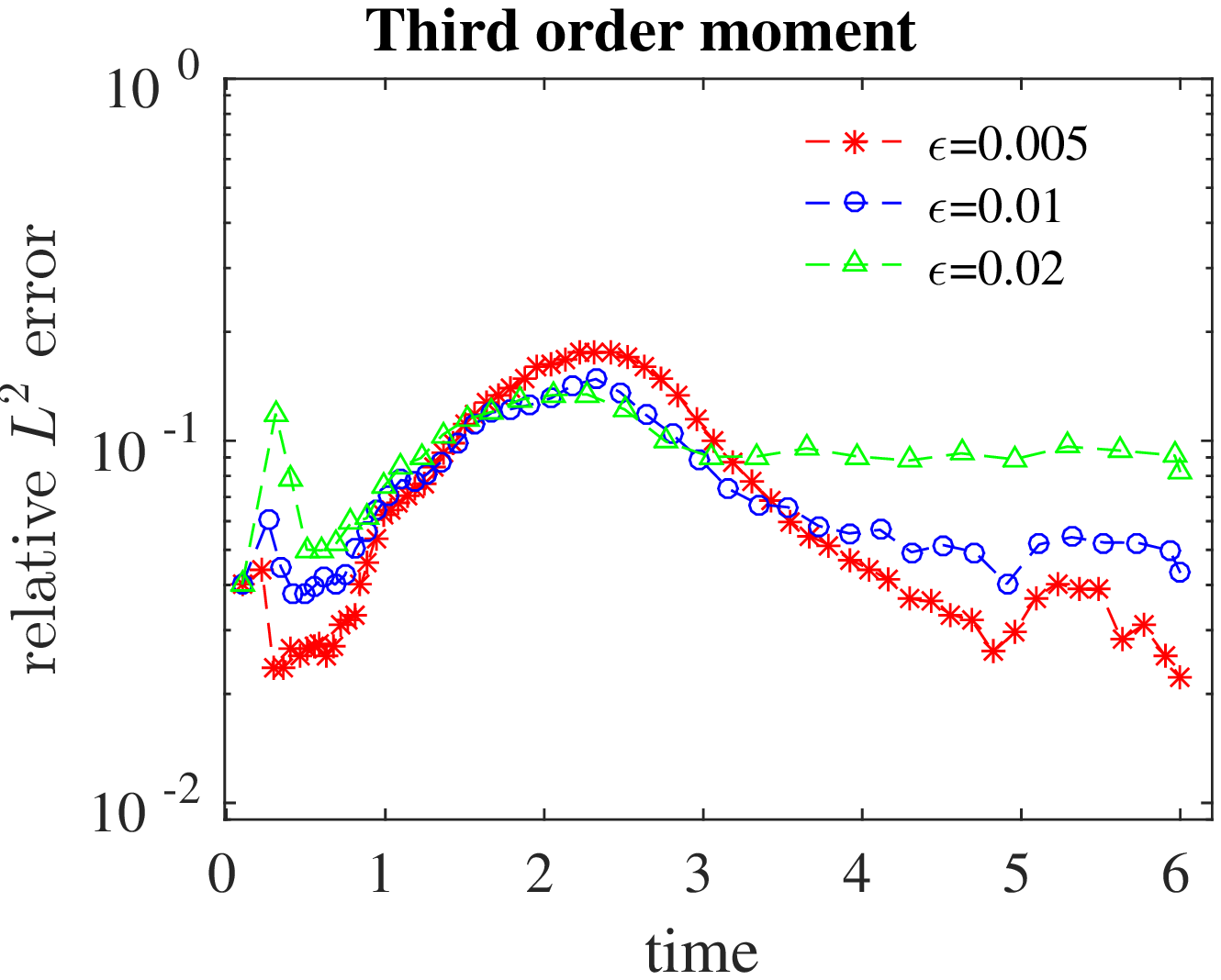}} 
\subfloat[]{
\includegraphics[scale=0.275]{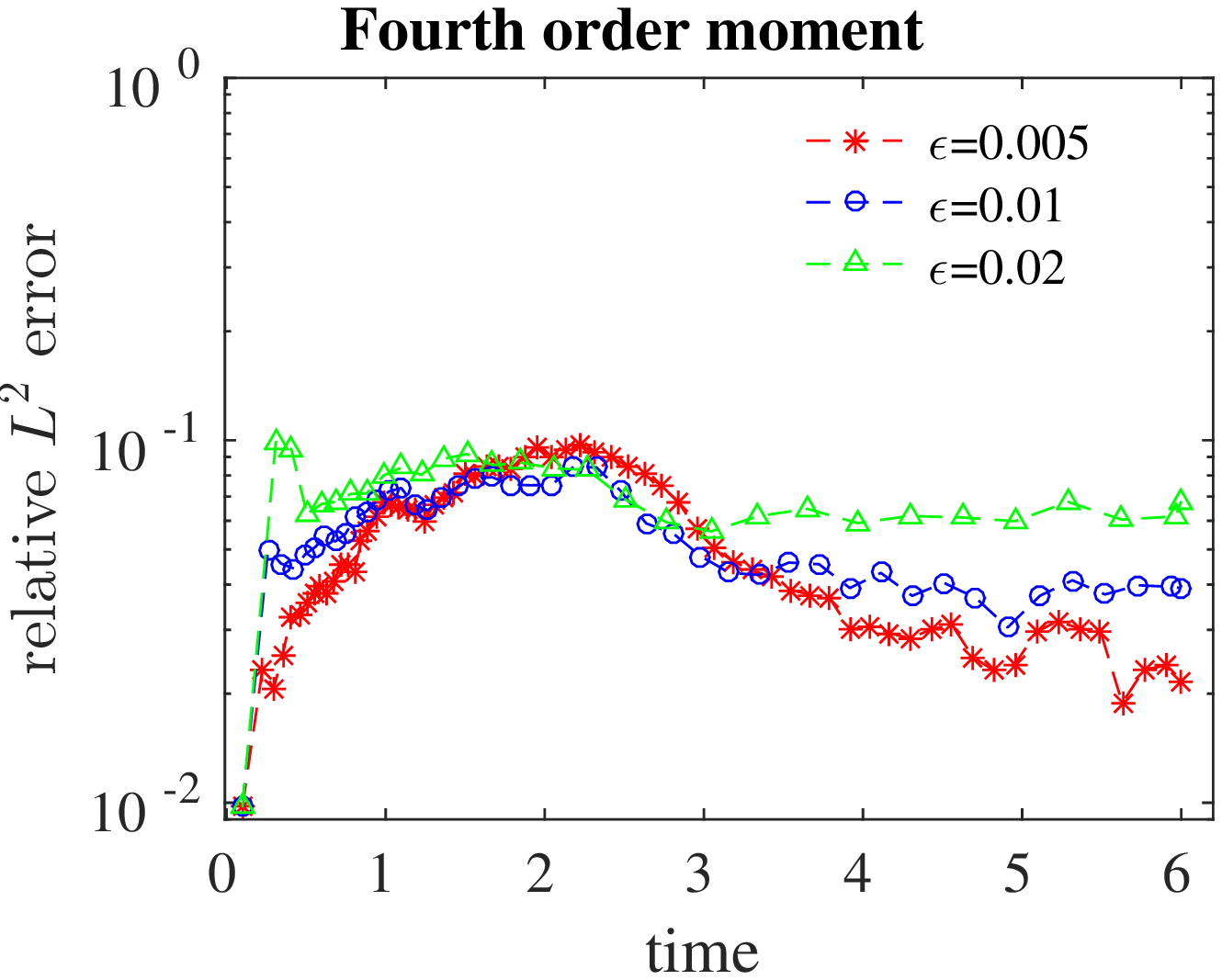}}  
\caption{Evolution of relative errors for moments with adaptive time stepping using different threshold values $\epsilon$. }
\label{fig:Ex1_3}
\end{figure} 

Using $\epsilon = 0.005, 0.01$ and $0.02$ results in $66,44$ and $29$ number of restarts with the time ratios $0.024$, $0.018$, and $0.015$ compared to MC3, respectively. As expected, decreasing the threshold value implies longer computational time and a larger number of restarts. Furthermore, from Figure \ref{fig:Ex1_3}, we observe that errors corresponding to $\epsilon=0.02$ are the largest in the initial stages and the long-term while errors for the smallest value $\epsilon=0.005$ correspond the most accurate behavior in the long-term. 

We make the following remarks: (i) optimal values of $\epsilon$ should be chosen according to the computational time and error level, and for this calculation, $\epsilon \in (0.005,0.01)$ seems optimal; (ii) earlier stages of the evolution should be analyzed carefully since using large values of the threshold value may result in a loss of accuracy; (iii) using very small values of $\epsilon$ may result in accumulation of errors if  errors at each restart are significant, e.g., when a small number of degrees of freedom is used in the KL expansion. Finally, we note that fitting a linear polynomial for \eqref{eq:rho} yields similar results. 

To show that the algorithm captures the invariant measure for a long time $T=6$, we consider the following three different initial conditions
\[
u_0(x) = 0.15 \sin( 2\pi x)  \quad \& \quad u_0(x) = 0.5 \cos( 4\pi x) \quad \& \quad u_0(x) = 0.25 (\sin(4\pi x)+\cos(8 \pi x)),
\]
and compute the moments at time $T=6$. We see from Figure \ref{fig:Ex1_2} that the dynamics converge to the unique invariant measure, which is a global attractor. 

\captionsetup[subfigure]{labelformat=empty}
\begin{figure}[!htb]
\centering
\subfloat[]{
\includegraphics[scale=0.275]{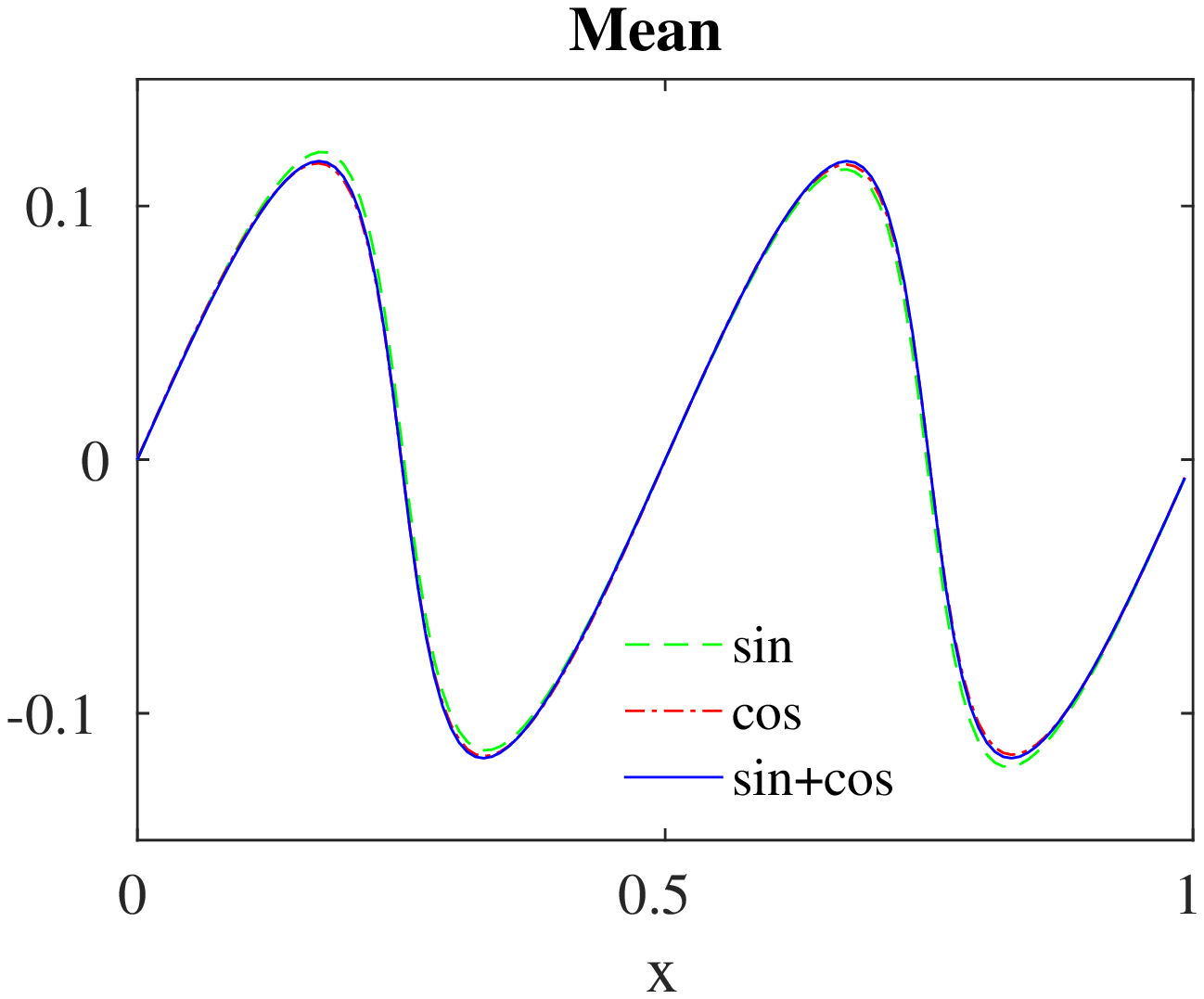}} 
\subfloat[]{
\includegraphics[scale=0.275]{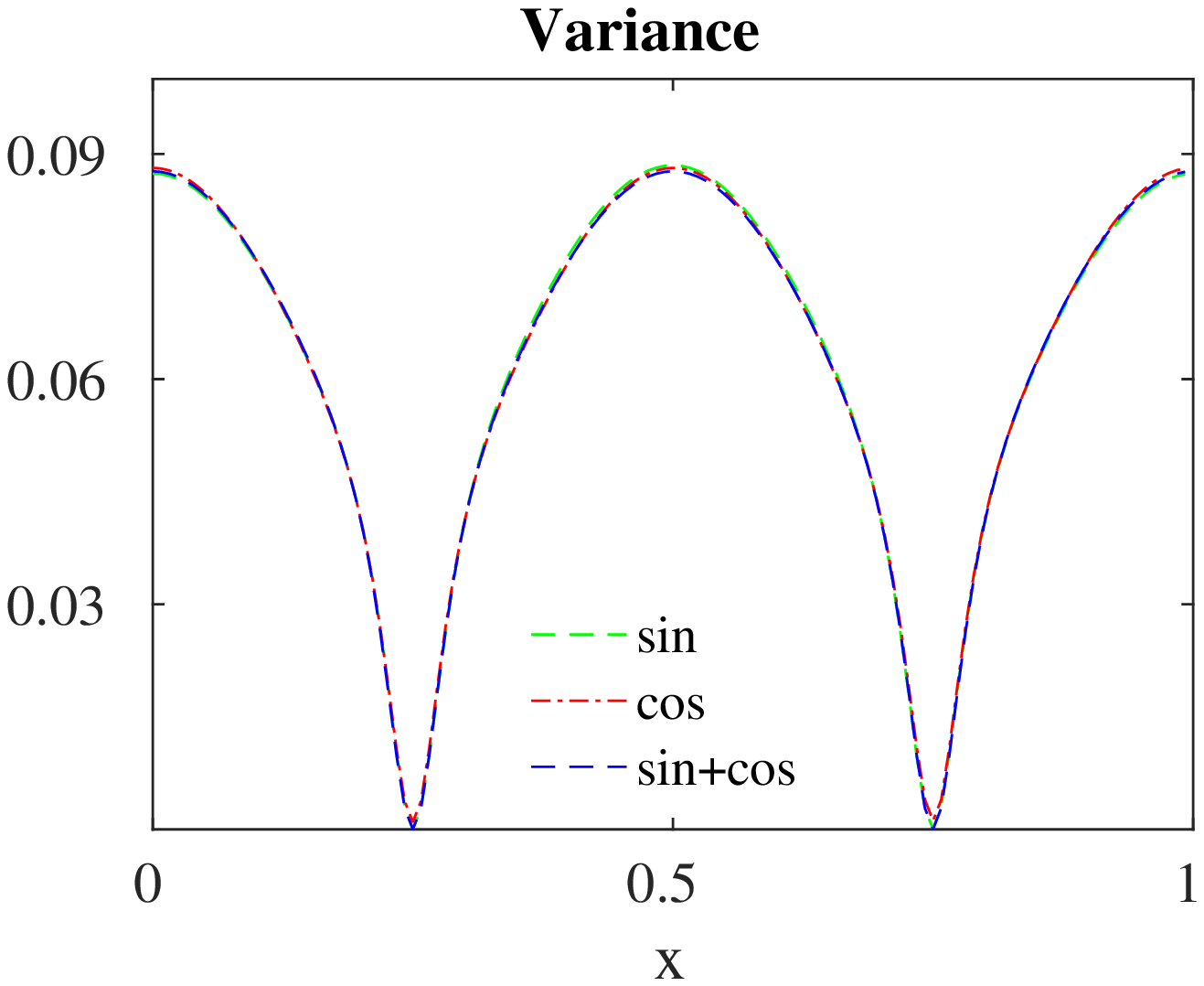}}  
\subfloat[]{
\includegraphics[scale=0.275]{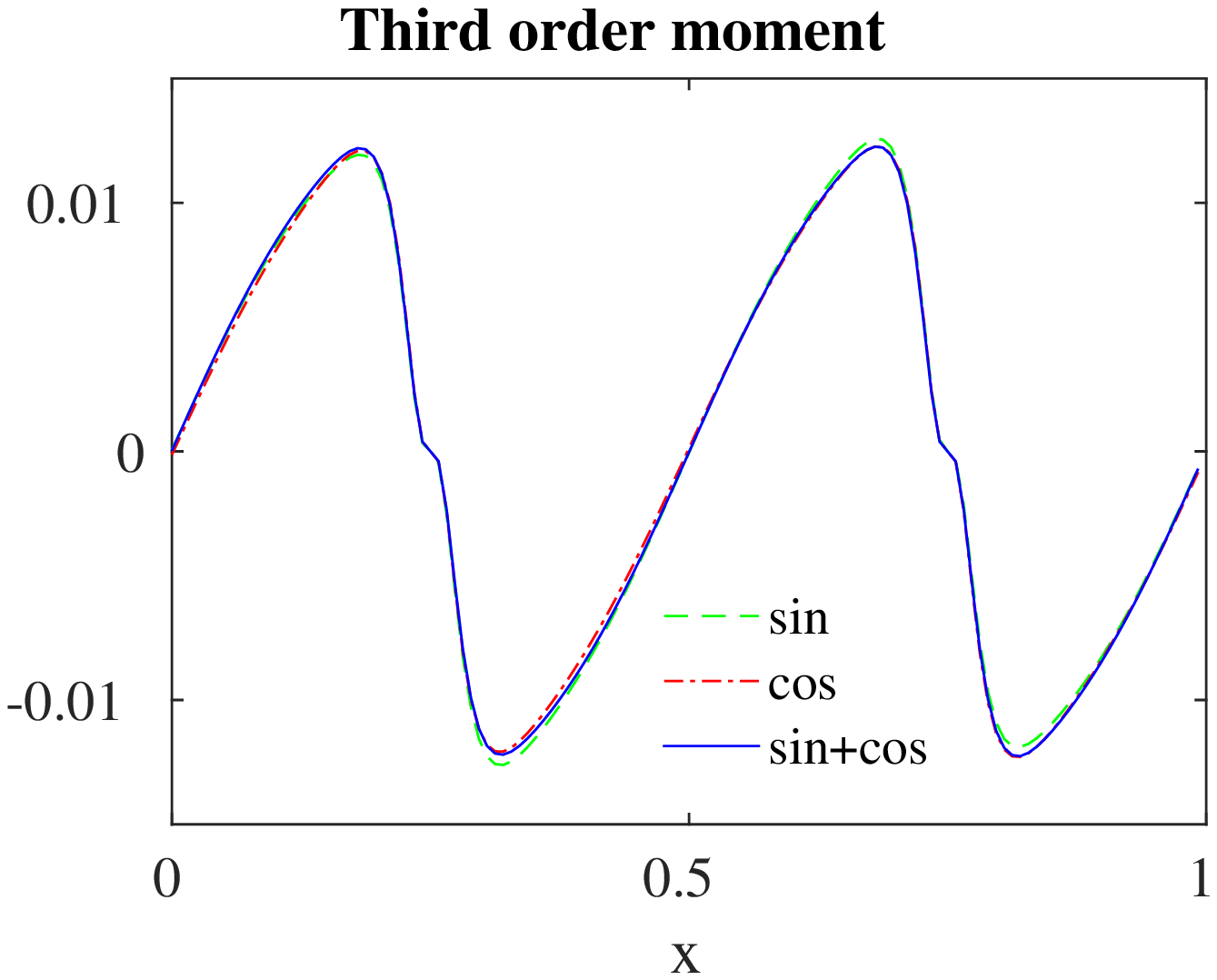}} 
\subfloat[]{
\includegraphics[scale=0.275]{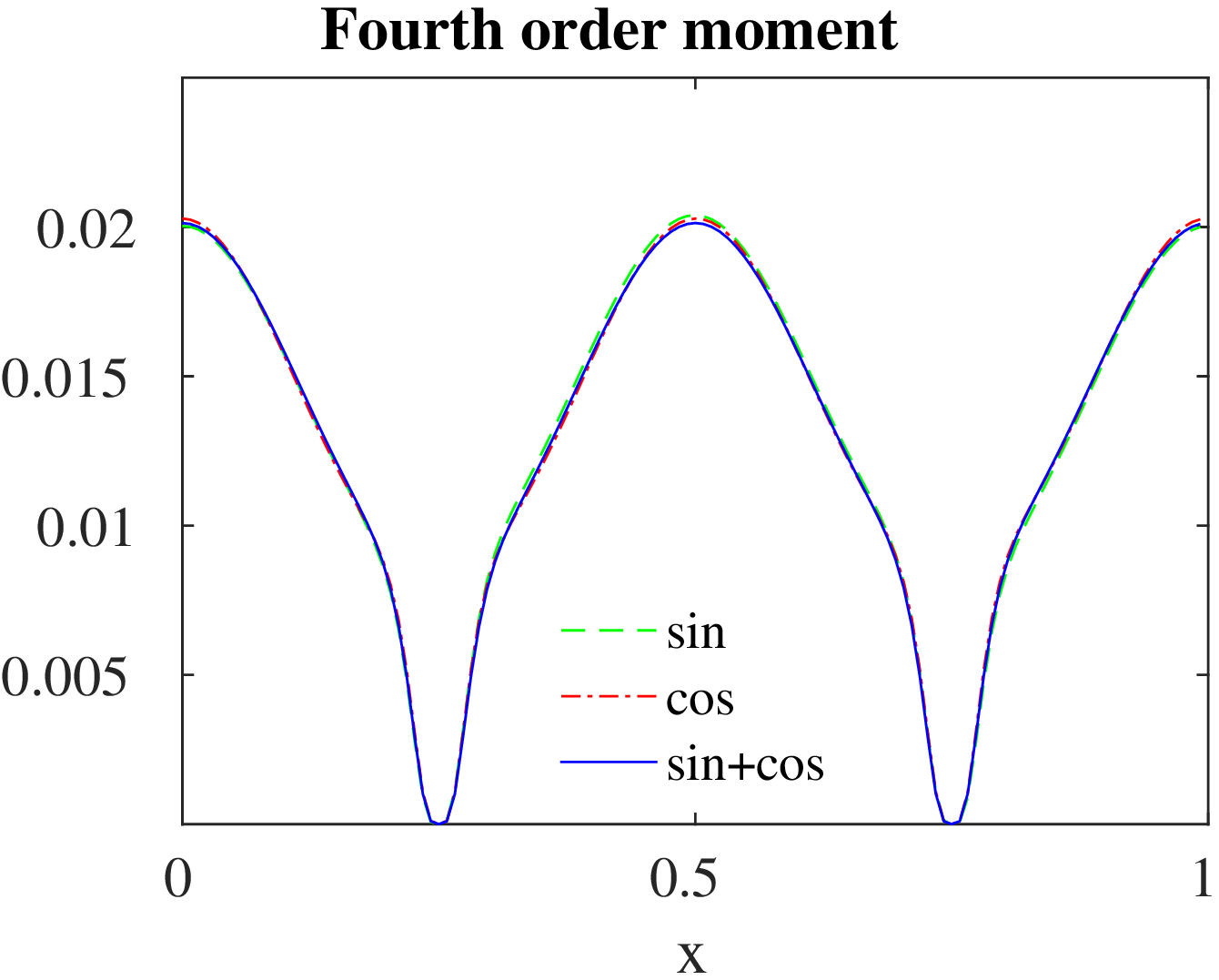}}  
\caption{Moments of the invariant measure at $T=6$ by DgPC. }
\label{fig:Ex1_2}
\end{figure} 
\end{exmp}

\begin{exmp} \rm
The purpose of the following numerical verification is to display the rate of convergence as $\Delta t$ varies for long time simulations. To this end, we take the initial condition and the forcing 
\[
u_0(x) = 0.5 (\exp( \cos (2 \pi x)  ) -1.5 )  \sin (2 \pi (x+0.37)) \quad \& \quad \sigma(x) = 0.5 \cos (4 \pi x),
\]
where the initial condition has several nonzero frequency components in the Fourier space. The viscosity is set to be $\nu= 0.005$. Note that there is no stationary state in the long-term. 

Using the same setting of Example \ref{ex:Burgers1} for parameters of DgPC, we apply DgPC with the varying values $\Delta t= 0.1, 0.2$ and $0.4 $ for each final time $T=2.4, 4.8, 9.6$ and  $14.4$. All simulations use the same time-step $dt= 0.001$ for time-integration.  Figure \ref{fig:Ex2_1} demonstrates that the order of convergence in long-time is varying between $O(\Delta t^{0.4})$ and $O(\sqrt{\Delta t})$. This behavior is consistent with the claims made in \cite[Theorem 5.1]{HLRZ06} in the setting of the Hermite PCE.

\begin{figure}[!htb]
\centering
\subfloat[$D=4$]{
\includegraphics[scale=0.3]{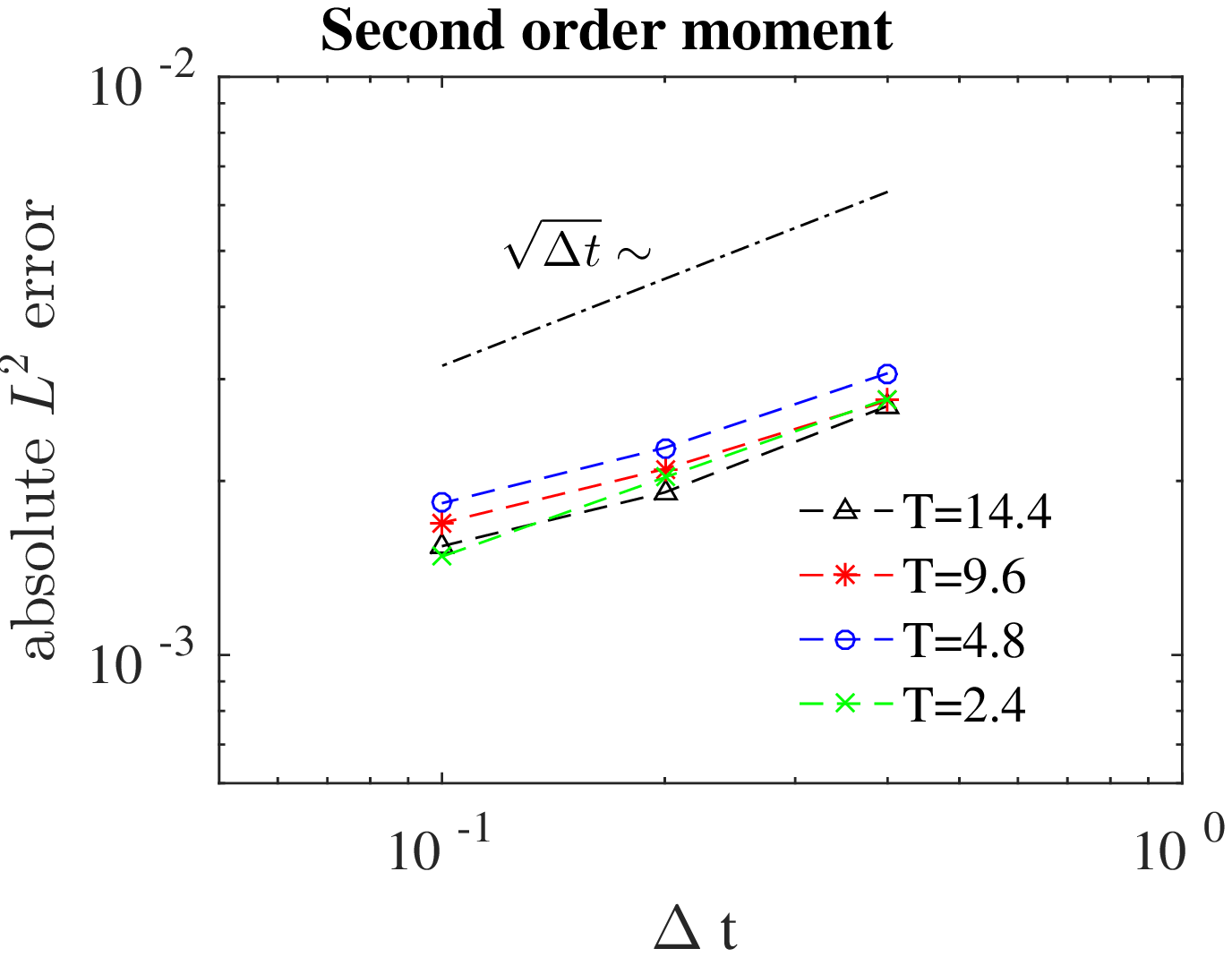}} 
\subfloat[$D=5$]{
\includegraphics[scale=0.3]{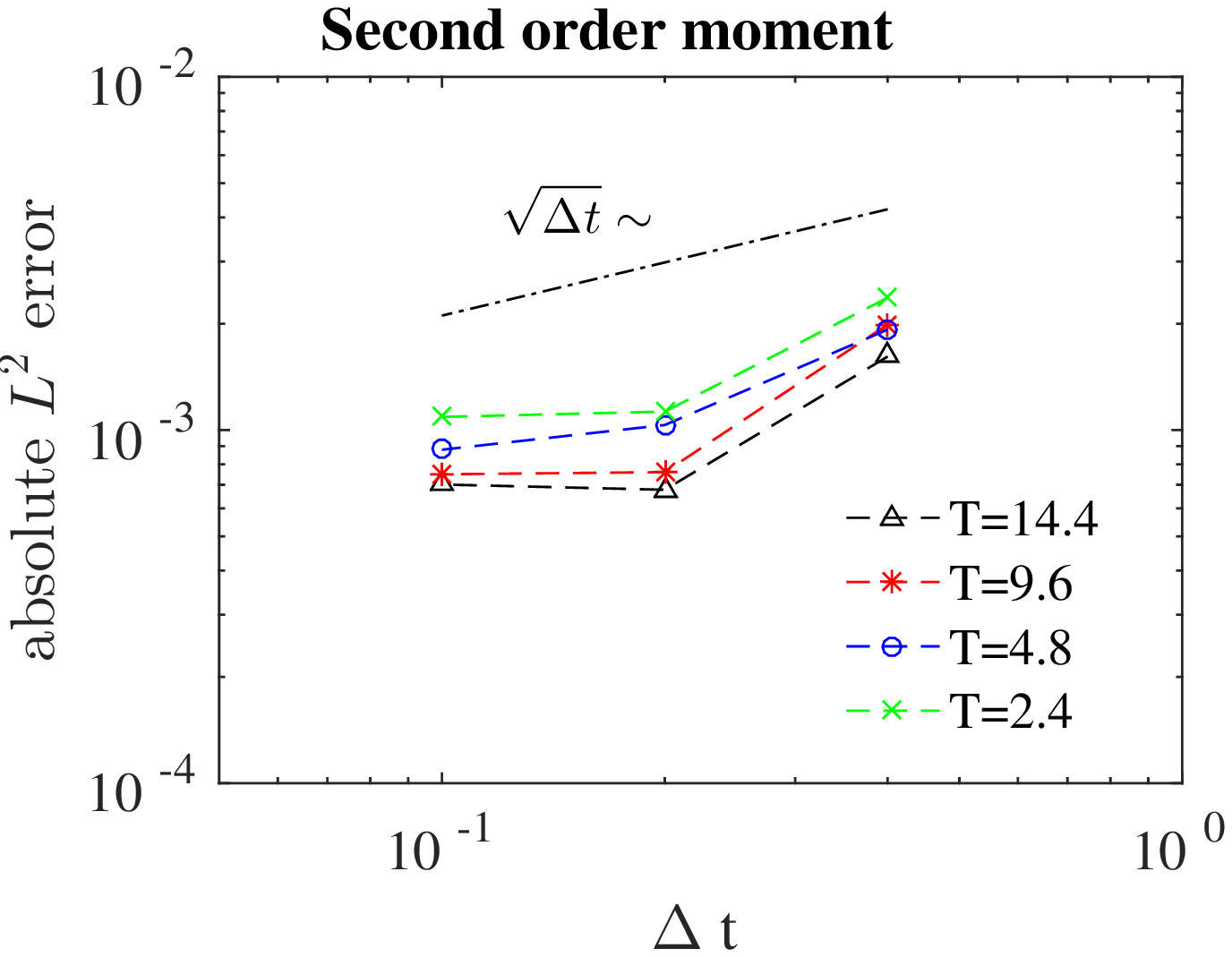}}  
\caption{Convergence behavior of errors in second moments using $\Delta t=0.1$, $0.2$, $0.4$. }
\label{fig:Ex2_1}
\end{figure}

\end{exmp}

\begin{exmp} \rm
In this example, we test the accuracy of the algorithm against an exact solution. When $\sigma(x) = \sigma$ is constant, the exact moments of the stochastic Burgers equation can be computed by solving the deterministic Burgers equation and estimating appropriate integrals by numerical quadratures.

We set $\sigma =0.1$ in \eqref{eq:1DBurgers} and take the initial condition 
\[
u_0(x) = 0.1-4  \nu \pi \cos(2 \pi x) /(3+\sin(2 \pi x)).
\]
In this case, the exact solution for the deterministic Burgers equation becomes 
\[  
u_{det}(x,t) = 0.1 - 4  \nu \pi \exp(-4 \nu \pi^2 t) \cos(2 \pi (x-0.1t))/ (3+ \exp(-4 \nu \pi^2 t) \sin(2 \pi (x-0.1t))), \quad t \geq 0. 
 \]
The moments of the stochastic solution $u(x,t)$ are then computed by the following integrals: 
\[
\E [u(x,t)^n] = \int_{\R^2} (u_{det}(x-z,t)+y)^n p(y,z) \, dy dz,
\]
where $p(y,z)= \frac{\sqrt{3}}{\pi \sigma^2 t^2}\exp \left( - \frac{2y^2}{\sigma^2 t} + \frac{6yz}{\sigma^2 t^2} - \frac{6z^2}{\sigma^2 t^3} \right) $ \cite[equation (3.13)]{HLRZ06}. We use a large number of quadrature points and the periodicity to compute the moments of the solution accurately.

To perform convergence analysis in terms of degree of polynomials, we take $N=1,2$ and $3$, and set $K=3$, $D=4$, $S=3 \times 10^5$. This setting results in $8,18$ and $38$ number of terms in the expansion for each time interval. Figure \ref{fig:ExConv_1} demonstrates the resulting relative errors for the moments of the solution. As expected, we observe that increasing the degrees of freedom $N$ helps to reduce the errors and an error accuracy of $O(10^{-3})$ is attained. Since $\sigma=0.1$ is held constant, the forcing term continuously forces the zeroth order spatial modes of  the high statistical moments, which are not damped by the viscosity and grow with time.

\begin{figure}[!htb]
\centering
\subfloat{
\includegraphics[scale=0.275]{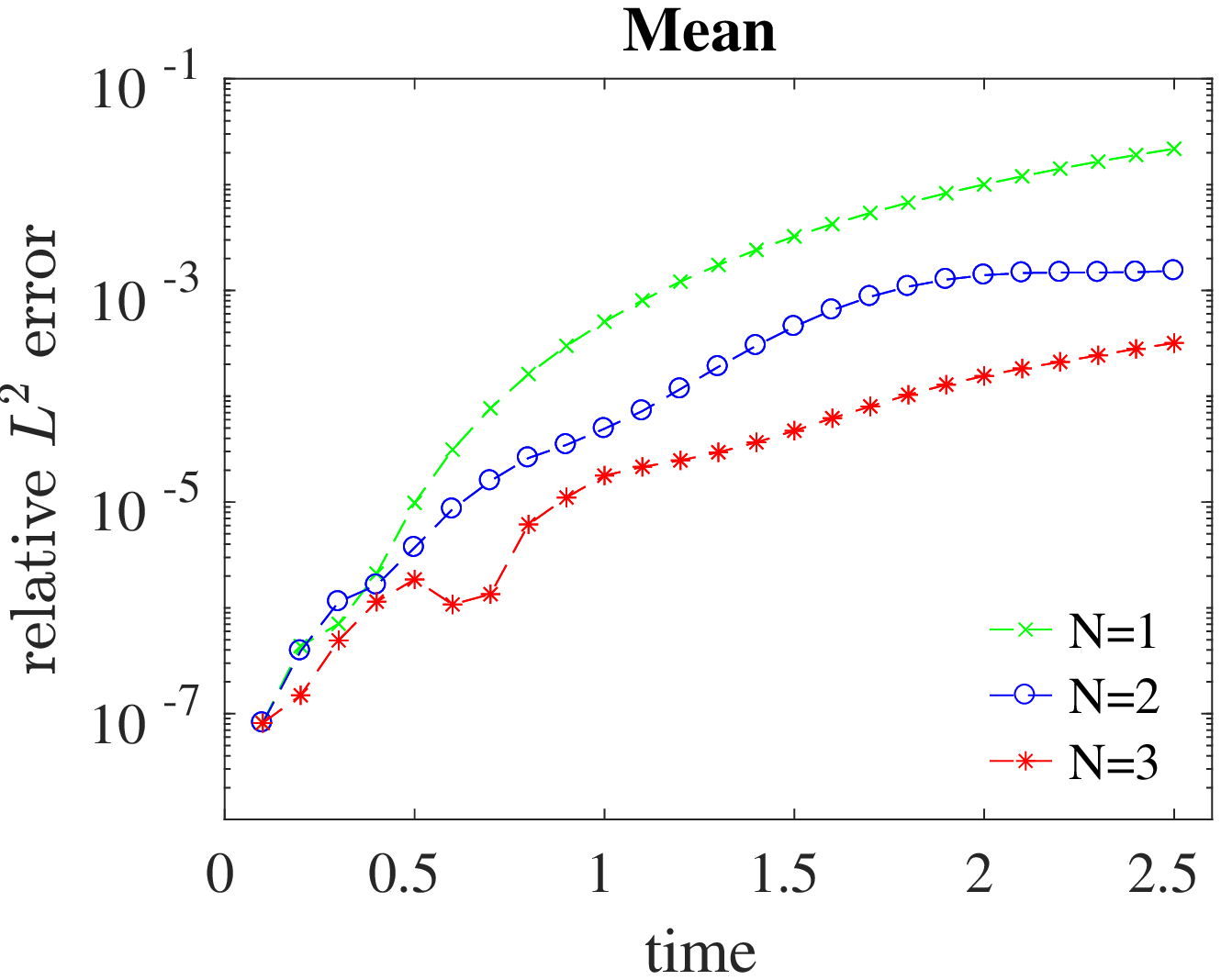}} 
\subfloat{
\includegraphics[scale=0.275]{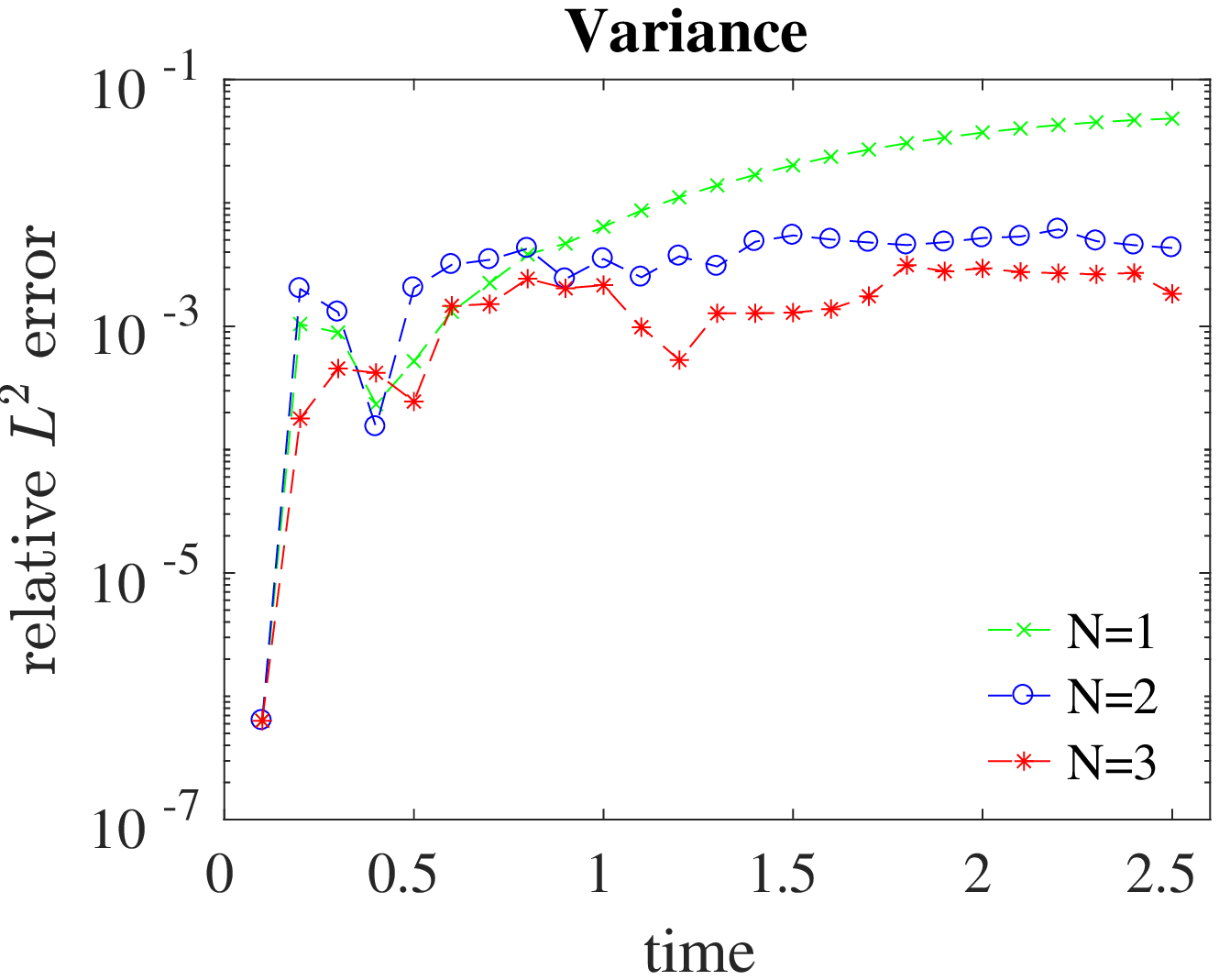}}  
\subfloat{
\includegraphics[scale=0.275]{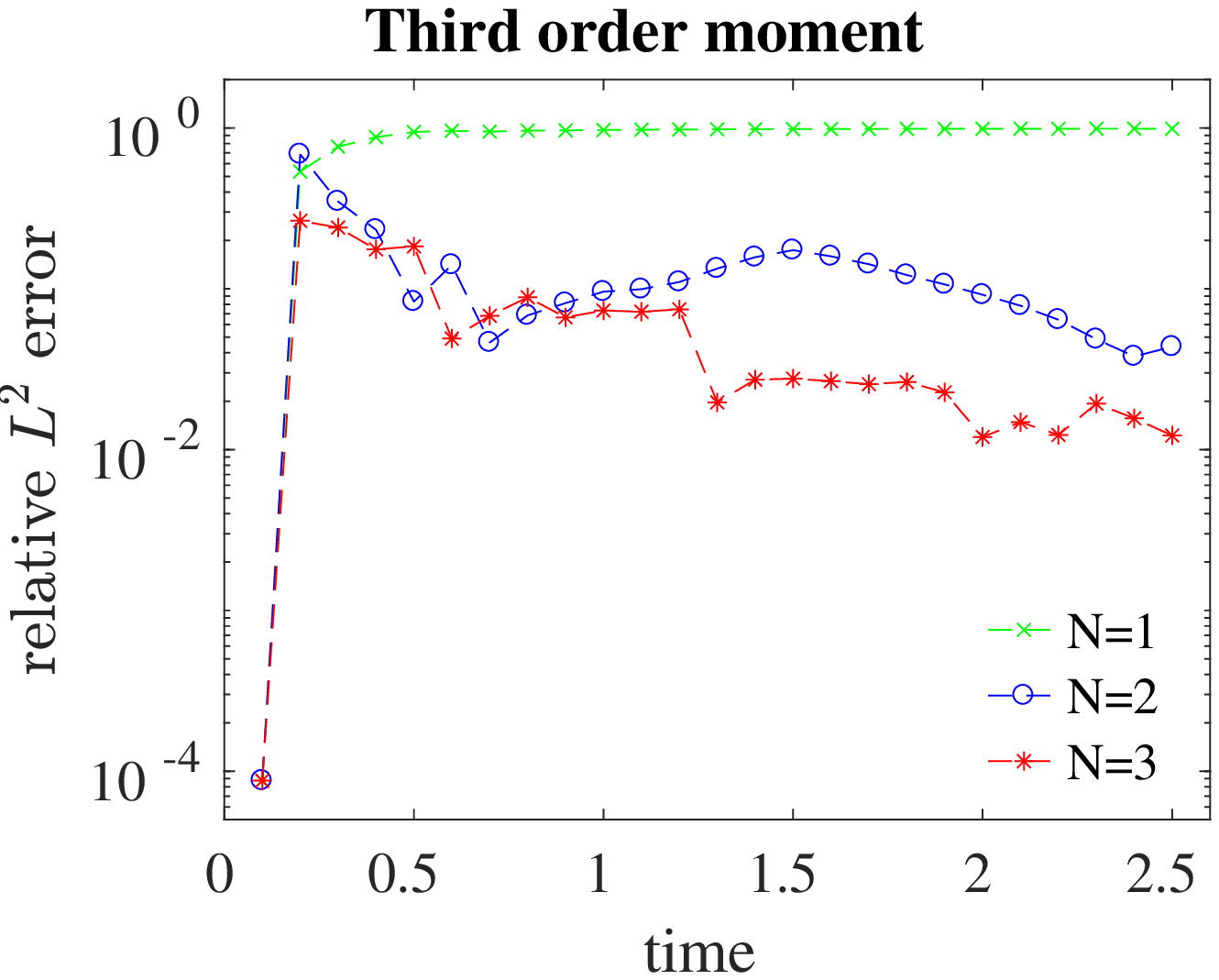}}
\subfloat{
\includegraphics[scale=0.275]{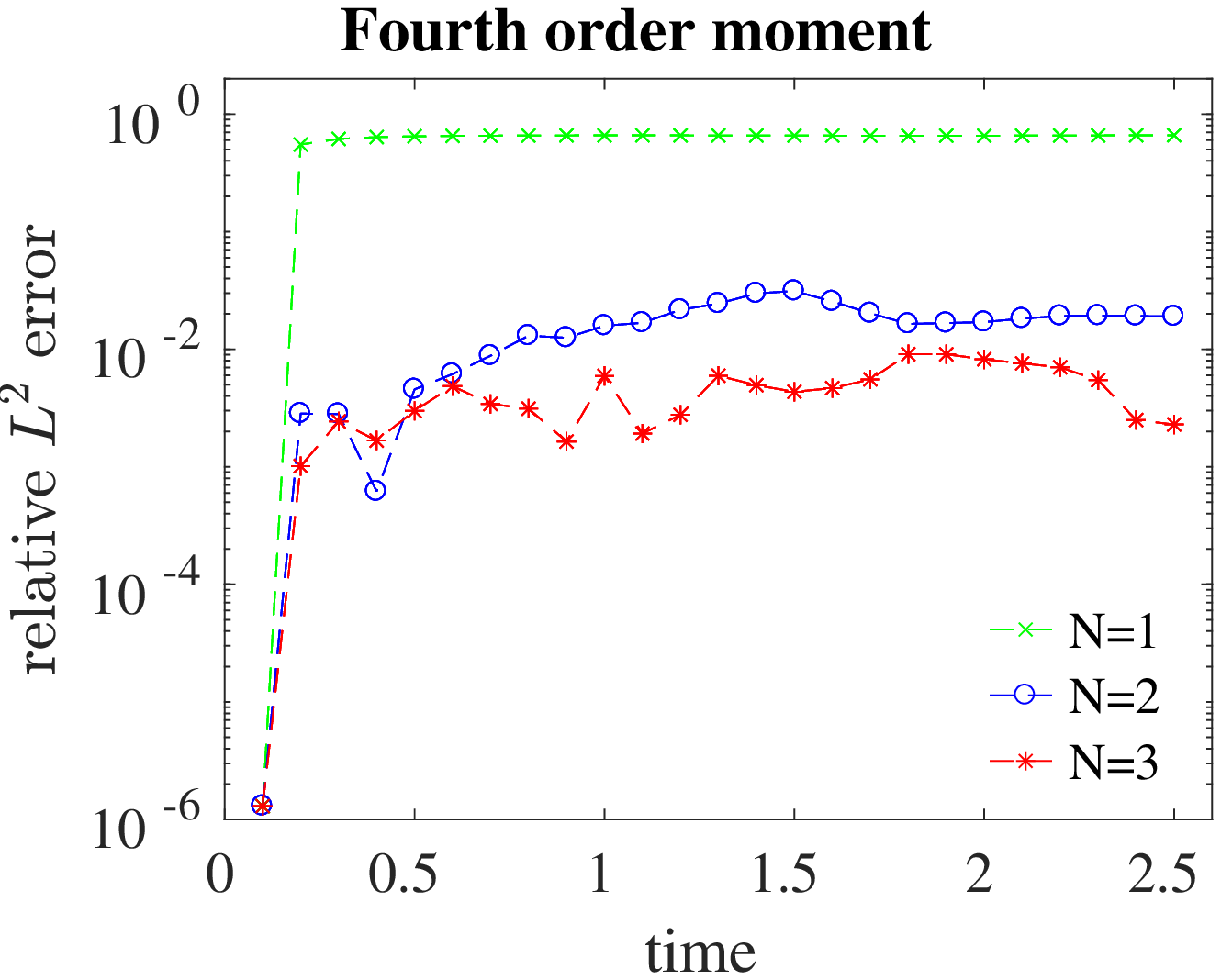}}    
\caption{Convergence behavior of errors in moments using polynomial degrees $N=1$, $2$, $3$. }
\label{fig:ExConv_1}
\end{figure}

\end{exmp}

\begin{exmp} \rm\label{ex:burg_het_vis}

For this simulation, we model the viscosity as an uncertain parameter with a spatial dependence, which can be useful for quantifying uncertainties in applications \cite{GS91, LeMK10,N09}. To this end, using the same notation of section \ref{sec:random_viscosity}, the covariance kernel of the underlying random process $Z(x,\omega)$ is assumed to be the following periodic exponential kernel
\[
\mbox{Cov}_{Z}(x,y) =   \sigma_Z^2 \exp \left( - \frac{2}{l_Z^2} {\sin^2(\pi (x-y))} \right), \quad x,y \in [0,1],
\]
where $\sigma_Z$ is the amplitude and $l_Z$ is the correlation length. We then compute the truncated KL decomposition of the mean zero process $Z(x,\omega)$
and  construct the viscosity as a function of $Z$ as
\[ 
\nu(x,\omega) = a_1 + Z^2(x,\omega), \quad \mbox{ where } \quad Z(x,\omega) = \sum_{l=1}^{D_{Z}} \sqrt{\lambda_l} \, U_l(\omega)  \, \phi_l(x),
\]
with independent uniform random variables $U_l \sim U(-1,1)$ and $a_1>0$. Armed with this viscosity, we consider the diffusion term in the Burgers equation \eqref{eq:1DBurgers} in divergence form, i.e., $ \partial_{x} (\nu(x,\omega) \, \partial_x u) $. 

We utilize Hermite and Legendre polynomial bases on the first subinterval $[0,t_1]$ to expand both Brownian motion and the random viscosity. Although the viscosity does not change in time, its PC representation changes as the PC bases differ at each restart time.  Thus, we keep track of the PCEs of both $u(x,t,\omega)$ and $\nu(x,\omega)$. 
Moreover, we compress both the solution and the random viscosity at each restart, i.e., the KL expansion is applied to the couple $(u_j,Z)$; see section \ref{sec:random_viscosity}. A reference calculation is computed by keeping the uniform random variables $U_l$, $l=1, \ldots ,D_Z$, at each restart in a PCE together with $\bxi_j$ and $\bta_j$. Relative errors are computed with respect to this reference calculation.

For the following simulation, the parameters are as follows: $K=2$, $ N=2$, $D=8$, $D_Z=3$, $S=3 \times 10^5$. The correlation length of the periodic kernel is set to $l_Z=2$. To avoid confusion, we slightly change the notation and denote by $\sigma_W$ the spatial part of the random forcing and consider three different scenarios:
\begin{enumerate}[ i)]
\item $\sigma_Z= 0.04$ and $\sigma_W(x) = 0.1\cos(2 \pi x)$, \label{sc:1}
\item $\sigma_Z= 0.1$ and $\sigma_W(x) = 0.1\cos(2 \pi x)$, \label{sc:2}
\item $\sigma_Z= 0.1$ and $\sigma_W(x) = 0.04\cos(2 \pi x)$,
\label{sc:3}
\end{enumerate}
with the same initial condition $u_0(x)=0.5\cos(4 \pi x)$.  These parameters correspond to different relative influences between the viscosity and the random forcing. 

We present the evolution of the relative errors for moments of the solution $u$ and the random viscosity $\nu$ for $T=4$ in  Figure \ref{fig:Ex3_1}. Each moment is averaged over distributions of the random process $Z$ and Brownian motion automatically by the PCE. First, we observe from the first two subfigures that the relative errors of the moments of the solution are of $O(10^{-2})$ in the long-time while the most accurate ones correspond to scenario \ref{sc:1}. In the same scenario, we see from the rightmost subfigure that the accuracy in the variance of $\nu$ decreases and stabilizes in time, which substantiates the observation that the algorithm selects the important part of the moments of the viscosity while keeping the solution accurate. Nevertheless, if slight changes in the moments of the viscosity become significant and accuracy needs to be improved, the KL expansion can be carried out using correlation matrices rather than covariance matrices.

Regarding the computational time, we note that the DgPC with $K+D=10$ number of variables  requires almost one eighth of the run time of the reference calculation which  utilizes $K+D+D_Z = 13$ number of variables in each PCE. Therefore, in cases where random parameters in the equation are high dimensional, applying the KL expansion to combined random variables is advantageous in terms of speed given a prescribed accuracy.

\captionsetup[subfigure]{labelformat=empty}
\begin{figure}[!htb]
\centering
\subfloat[]{
\includegraphics[scale=0.275]{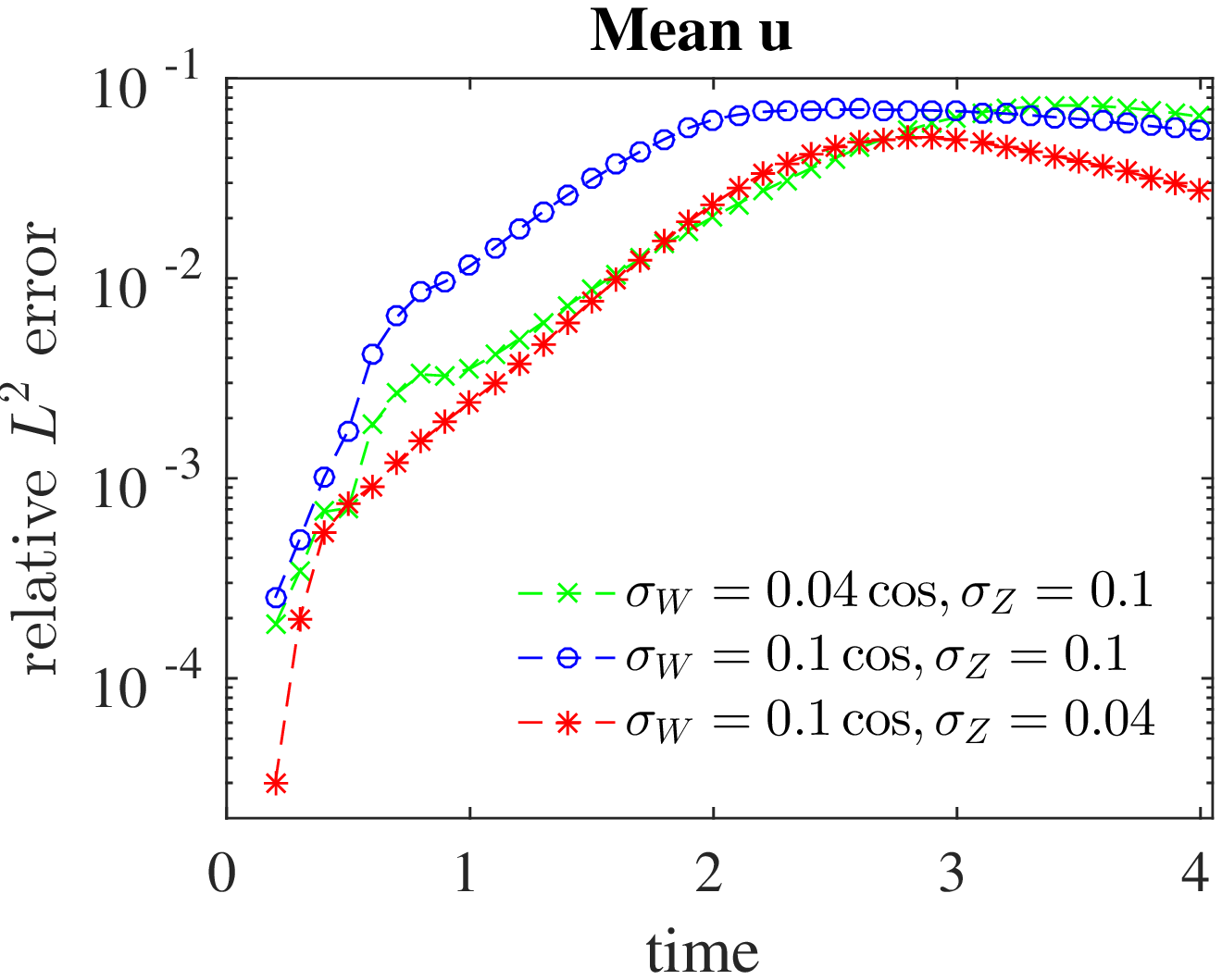}} 
\subfloat[]{
\includegraphics[scale=0.275]{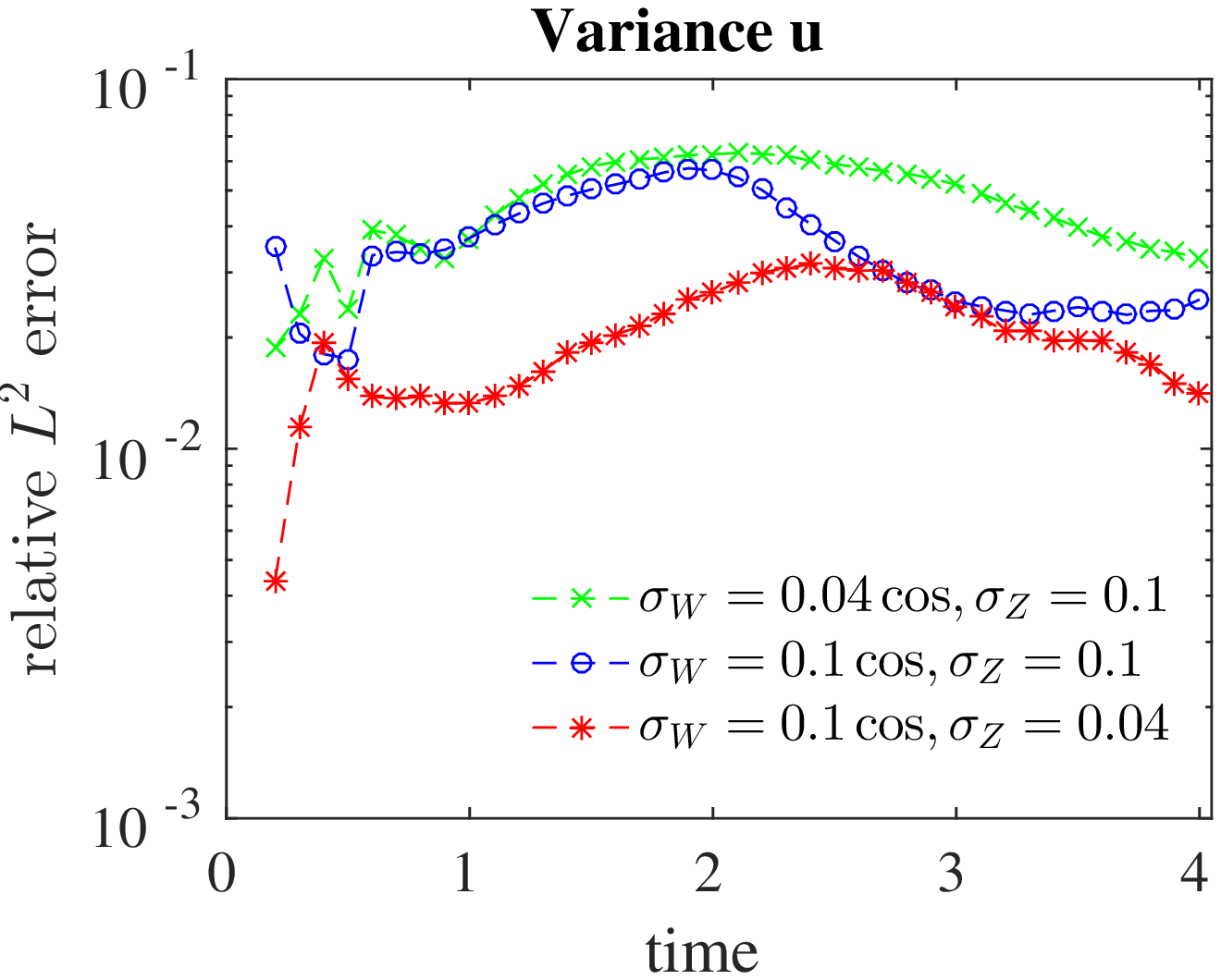}} 
\subfloat[]{
\includegraphics[scale=0.275]{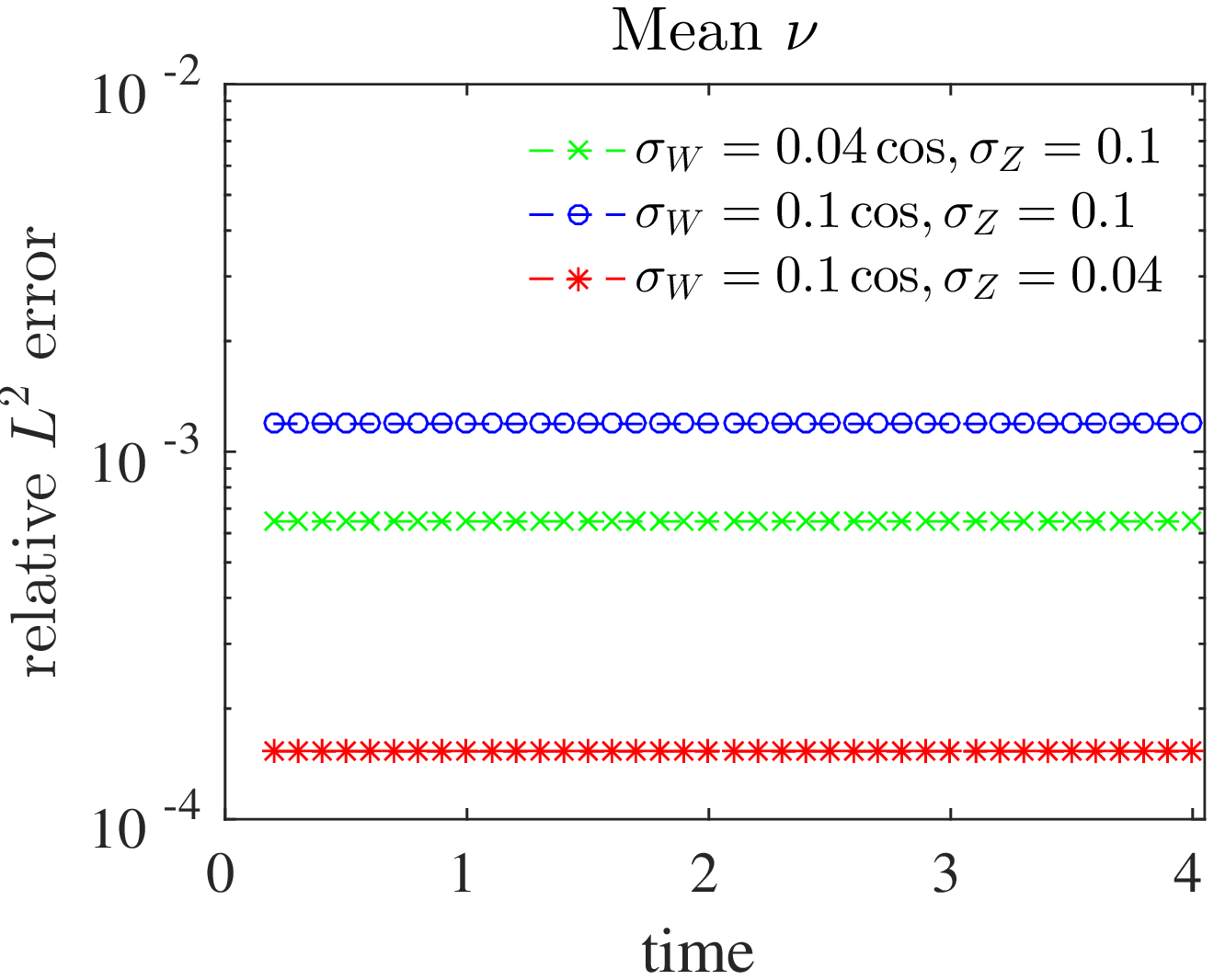}} 
\subfloat[]{
\includegraphics[scale=0.275]{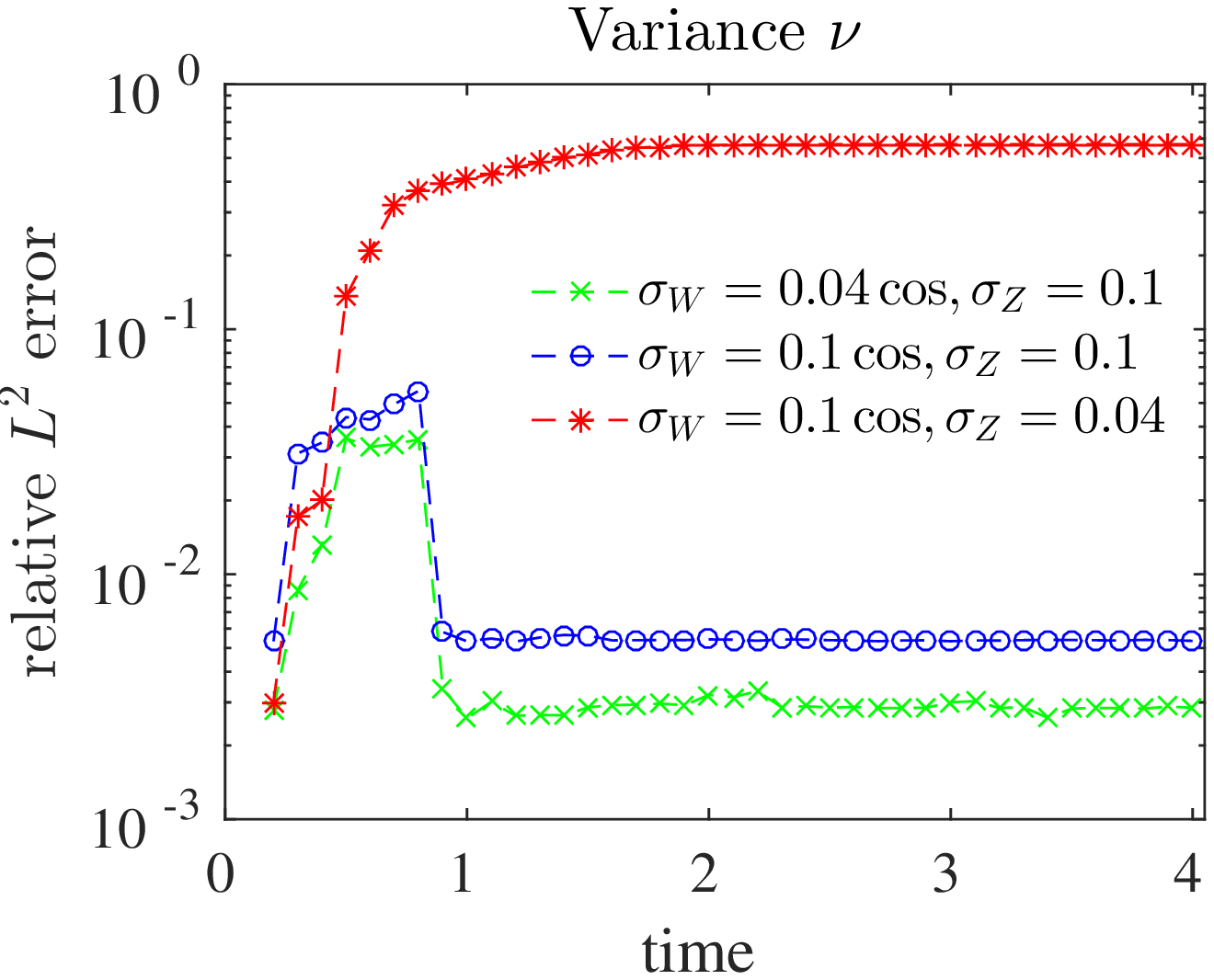}}  
\caption{Relative errors for moments of the solution and the random viscosity. Errors are computed by comparing DgPC with $D=8$ to a reference calculation which uses $D=8$ and $D_Z=3$. }
\label{fig:Ex3_1}
\end{figure} 

Figures \ref{fig:Ex3_2a} and \ref{fig:Ex3_2b} show forty snapshots of the moments of the solution in time for the scenario \ref{sc:3}, where black curves indicate the initial states for each moment. Note that since the viscosity is random, for each realization of the viscosity there is a unique invariant measure, and what these figures exhibit is the steady state, which is obtained by averaging these measures over the distribution of the viscosity. We also see that the moments of this averaged measure differ in magnitude compared to those in Figure \ref{fig:Ex1_2} as the relative influences of viscosity and Brownian motion are changed. The rightmost sub-figure plots the one-point probability density function corresponding to this steady state. The density function is easily obtained using the samples of the approximated random field via a kernel density estimation procedure; see Remark \ref{rmk:samples} and \cite{BGK10}. We see that for each point $x \in [0,1]$, the density function is unimodal and has peaks near the points $x$, where the variance is minimum. We finally note that the cross covariance structure of the solution on the spatial mesh can also be easily deduced from the algorithm if needed. 

\captionsetup[subfigure]{labelformat=parens}
\begin{figure}[!htb]
\centering
\subfloat[Mean over time]{\label{fig:Ex3_2a}
\includegraphics[scale=0.3]{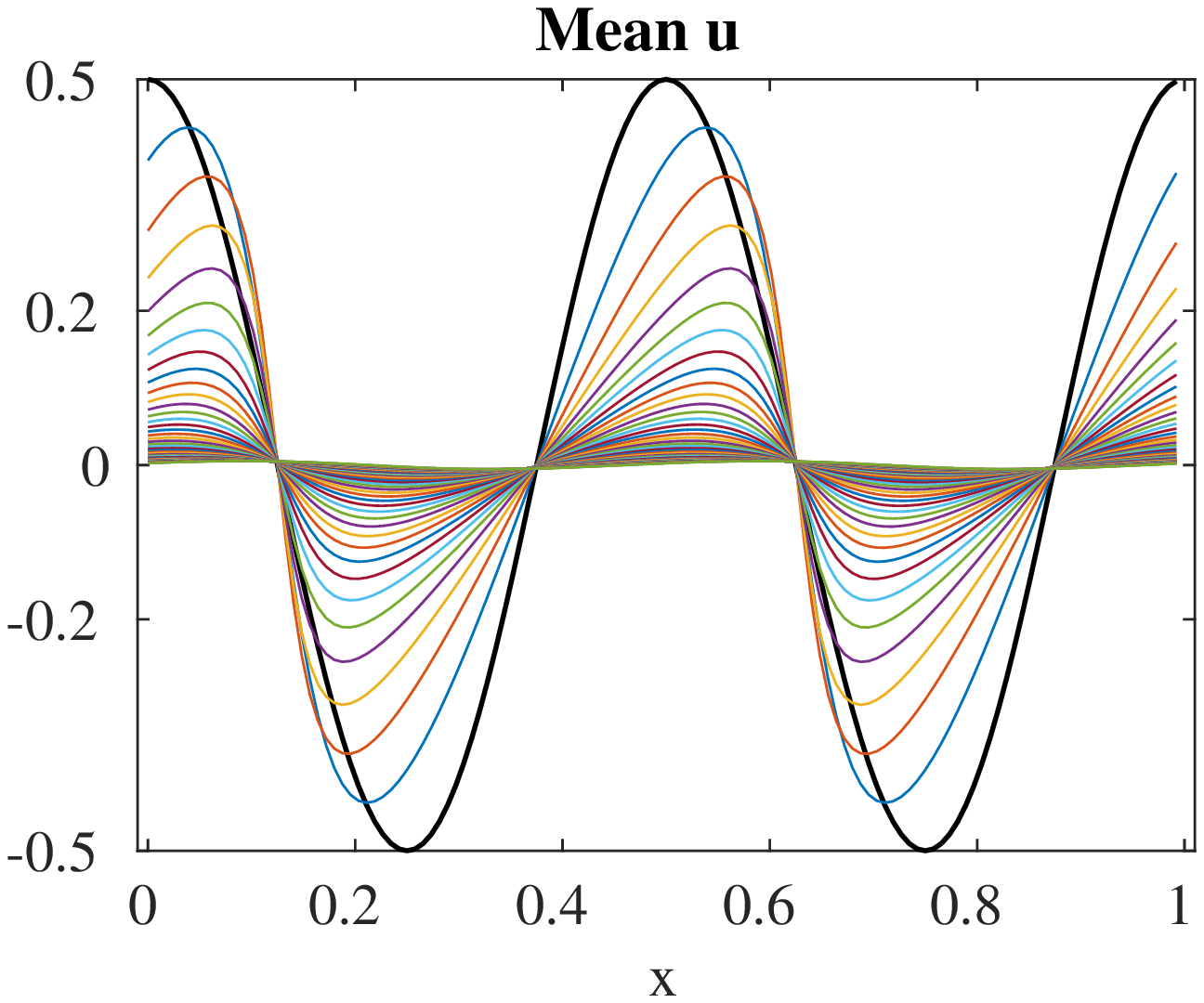}}   
\subfloat[Variance over time]{ \label{fig:Ex3_2b} 
\includegraphics[scale=0.3]{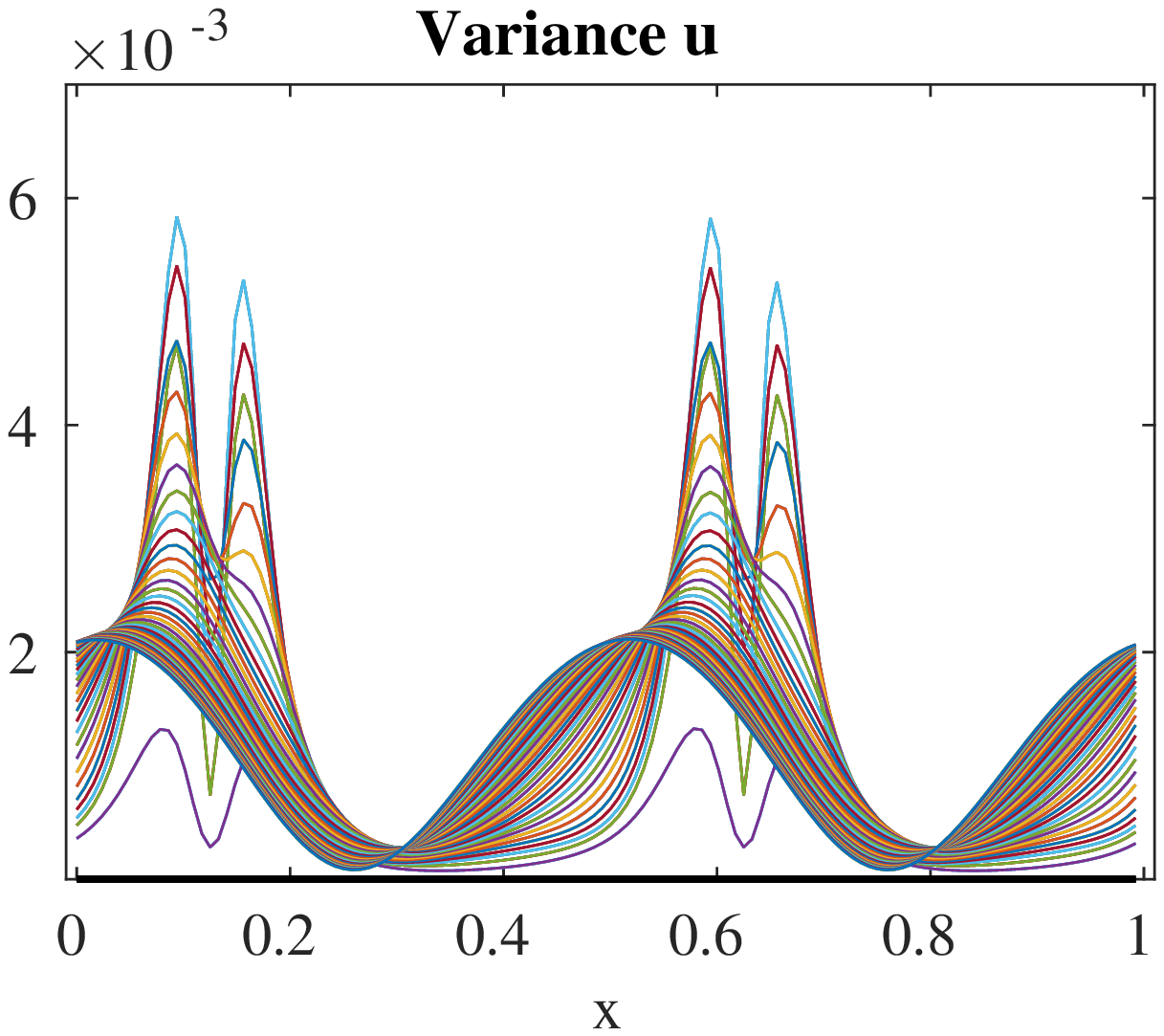}} 
\subfloat[One-point pdf at final time]{\label{fig:Ex3_2c}
\includegraphics[scale=0.3]{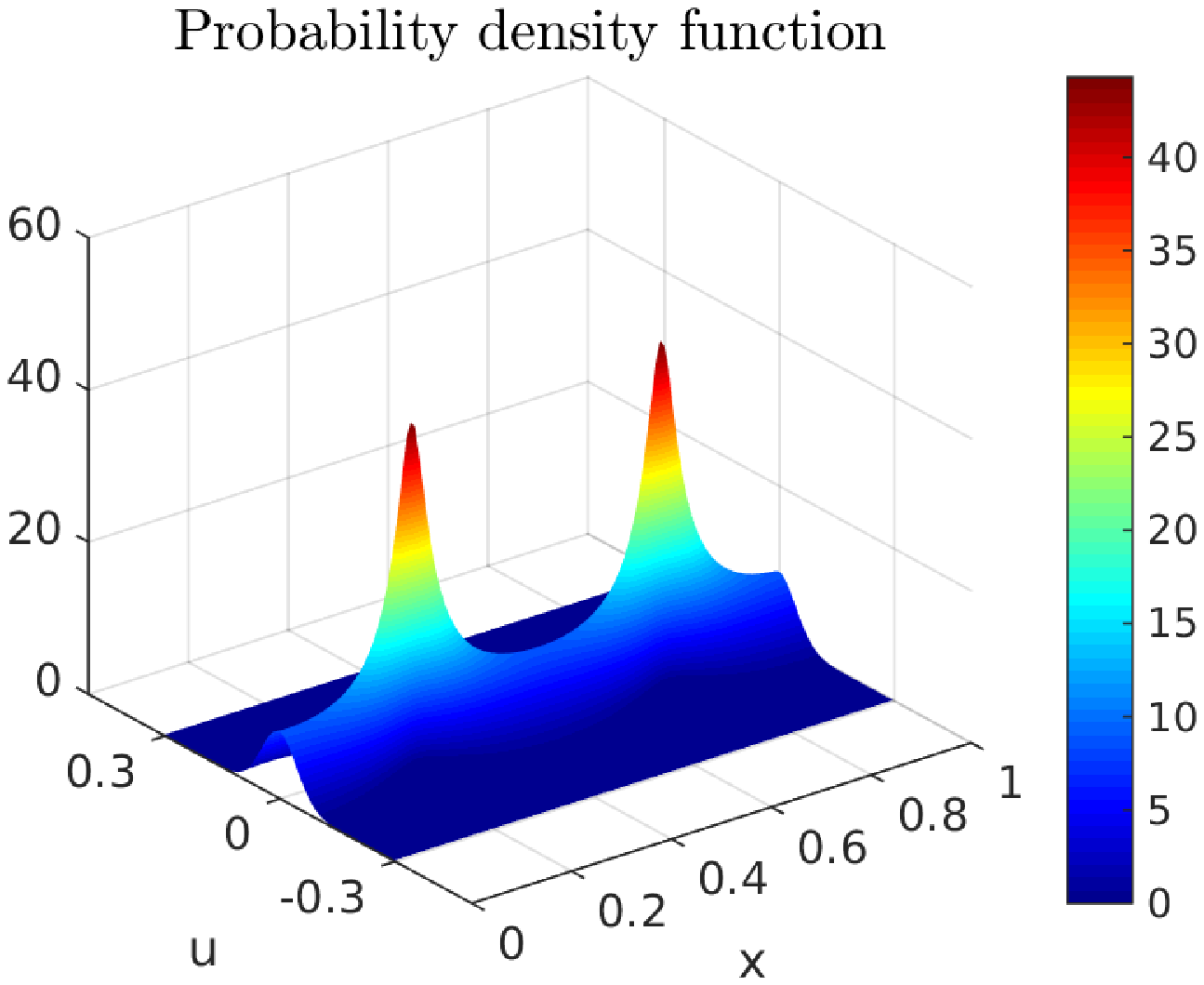}}  
\caption{Snapshots of second order moments on $[0,T]=[0,4]$ and one-point probability density function at $T=4$.}
\label{fig:Ex3_2}
\end{figure}

\end{exmp}

\subsection{Navier--Stokes equations} \label{sec:SNS}

In this section, we provide applications of our algorithm to solve a two-dimensional system of stochastic Navier--Stokes equations (SNS). We consider the following coupled system:
\begin{align} \label{eq:SNS}
\begin{cases}
 \theta_t +  \mathbf{u} \cdot \nabla \theta = \mu \, \Delta \theta, \\
 \mathbf{u}_t +\mathbf{u}\cdot \nabla \mathbf{u} = \nu \Delta \mathbf{u} - \nabla P + \sigma \dot{\mathbf{W}}(t), \\
 \nabla \cdot \mathbf{u} =0. 
 \end{cases}
\end{align}
where $\theta$, the temperature, is convected by the stochastic velocity field  $\mathbf{u} =(u,v)$, which is forced by a Brownian motion $\mathbf{W} =(W_1,W_2)$ with independent components and the spatial part $\sigma(x,y)= \mbox{diag}(\sigma_1(x,y),\sigma_2(x,y))$, $x,y \in [0,1]$. The temperature diffusivity and fluid viscosity are denoted by $\mu$ and $\nu$, respectively. Notice that equations are coupled only through velocity term and the temperature is convected by the random velocity passively.

We take the computational domain $[0,1] \times [0,1] \subset \R^2$, and  assume that $\theta$ and $\mathbf{u}$ are doubly periodic. It is then possible to introduce the stream function $\psi$ such that $\mathbf{u} = (\psi_y, - \psi_x)$, define the vorticity $w = v_x - u_y$ and rewrite the system \eqref{eq:SNS} as: 
\begin{align} \label{eq:sns_vor}
\begin{cases}
 \theta_t + (u \theta)_x + (v \theta)_y = \mu \, \Delta \theta, \\
w_t + (u w)_x + (v w)_y = \nu \Delta w + (\sigma_2)_x \dot{W_2}- (\sigma_1)_y \dot{W_1} , \\
-\Delta \psi = w, \quad 
u = \psi_y, \quad  v= - \psi_x.
 \end{cases}
\end{align}
We also suppose that the stream function is periodic. Following \cite{HLRZ06,Luo}, the initial condition for the vorticity $w$ is chosen to be 
\begin{align} \label{init_vor}
w(x,y,0) = C - \frac{1}{2 \delta} \exp  \left ( - \frac{I(x)(y-0.5)^2}{2 \delta^2} \right),
\end{align}
where $I(x) = 1+ \varepsilon (\cos (\gamma \, 2 \pi x) -1)$ with $\gamma \in \mathbb{N}$, and $C$ is a constant to make the initial condition mean zero on $[0,1]^2$. This choice corresponds to a flat shear layer of width $\delta$ concentrated at $y=0.5$. The width is  perturbed sinusoidally with the amplitude $\varepsilon$ and spatial frequency $\gamma$.  In our numerical experiments, we will also consider the reflected initial condition $w(y,x,0)$ for which the layer is concentrated vertically;  see Figures \ref{fig:Ex4_1} and \ref{fig:Ex7_1}.    

We consider the following initial condition for the temperature: 
\begin{align*} 
\theta(x,y,0) = \begin{cases}
 H_{\delta}(y-0.25),  & \mbox{if } y \leq 0.4, \\
 1- 2H_{\delta}(y-0.5), & \mbox{if }  0.4 < y < 0.6 ,\\
 -H_{\delta} (0.75-y), & \mbox{if } y \geq 0.6,
 \end{cases} \quad \mbox{where} \quad  H_{\delta}(x) = \begin{cases}
0, & \mbox{if } x < -\delta,\\
\frac{x+\delta}{2\delta} + \frac{\sin(\pi x /\delta)}{2\pi}, & \mbox{if }  |x| \leq \delta, \\
1, & \mbox{if } x > \delta.
 \end{cases}
\end{align*}
is the mollified Heaviside function. This formulation yields an initial condition, which consists of four connected layers, where interfaces between layers have thickness $\delta$. As discussed in \cite{HLRZ06,Luo}, setting small values to $\delta$ creates a sharp shear layer and temperature interface.  As a consequence of Kelvin-Helmholtz instability, the fluid then will roll-up, and the temperature will be convected and mixed up \cite{LSM93}.  

The Hermite PCE is effectively applied to stochastic Navier--Stokes equations (in most cases without Browninan motion forcing) in various manuscripts \cite{LKHG01, Xiu_thesis, XK03, WK05, HLRZ06, Luo}. The presence of Brownian motion forcing makes the system very hard to solve even for short times due to the overwhelming number of random variables needed in the Hermite PCE~\cite{HLRZ06,OB16}. Therefore, we now apply the DgPC algorithm to the system \eqref{eq:sns_vor}.

We choose the same orthonormal system $m_{j,i}(t)$ on $[t_j,t_{j+1}]$ as in the previous section and project each component $W_1(t)$ and $W_2(t)$ of the Brownian motion to obtain $\bxi_j$. The total number of $\xi_i$'s will be denoted by $K$, where the first (last) $K/2$  variables correspond to the first (second) component of $\mathbf{W}(t)$. Assuming the variables $w, \theta, \mathbf{u}$ and $\psi$ admit PCEs, the method keeps track of the corresponding expansions. At each time-step $t_j$, the KL decomposition is applied to $(w_j, \theta_j)$ and the mode $\bta_j$ is obtained. Then, Galerkin projection of \eqref{eq:sns_vor} onto the space $ \mbox{span} \{ T_{\balp}(\bxi_{j},\bta_j) : \,\alpha \in \mathcal{J}^r_{K+D,N} \}$ yields the following equations: 
\begin{align} \label{eq:sns_vor_sys}
\begin{cases}
\partial_t (\theta_{j+1,\balp}) +  \partial_x (u_{j+1} \theta_{j+1})_{\balp} +  \partial_y (v_{j+1} \theta_{j+1})_{\balp} = \mu \, \Delta \theta_{j+1,\balp}, \\
 \partial_t( w_{j+1,\balp}) + \partial_x (u_{j+1} w_{j+1})_{\balp} +\partial_y   (v_{j+1} w_{j+1})_{\balp} = \nu \, \Delta w_{j+1,\balp} \\  
 \qquad \quad \quad  +  (\sigma_2)_x \displaystyle  \sum_{i=K/2+1}^K m_{j,i-K/2}(t)\, \E[ \xi_{j,i} \, T_{\balp} ]   - (\sigma_1)_y \sum_{i=1}^{K/2} m_{j,i}(t) \, \E[ \xi_{j,i}\, T_{\balp} ] , \\
-\Delta \psi_{j+1,\balp} = w_{j+1,\balp}, \quad 
u_{j+1,\balp} =  \partial_y (\psi_{j+1,\balp}), \quad  v_{j+1,\balp}= - \partial_x ( \psi_{j+1,\balp}).
 \end{cases}
\end{align}
The resulting deterministic PDE system \eqref{eq:sns_vor_sys} is solved utilizing a truncated Fourier series and FFT in two dimensions on a mesh of size $M \times M$.

 As we discussed in Section \eqref{sec:KLexpansion}, there are two main methods to compute the KL expansion: (i) assemble the full covariance matrix using \eqref{eq:cov_pce} and use a Krylov subspace method to find largest eigenvalues (as was done in the previous section); or (ii) use the random projection technique to accelerate the computation by equation \eqref{eq:cov_QQ'} and find eigenvalues of the resulting small matrix. To show the computational savings incurred by the second method, some SNS systems were solved using both methods and accuracies are compared.

\begin{exmp} \label{ex:sns_short1} \rm 

This simulation concerns the short-time accuracy and the computational time of the DgPC algorithm, which will be assessed using comparisons with Monte Carlo methods with a sufficiently high number of samples.

 We set $\mu,\nu =0.0002$ and take the spatial part of forcing as
\[
(\sigma_1)_y = 0.1 \pi \cos (2 \pi x) \cos (2 \pi y), \quad (\sigma_2)_x = 0.1 \pi \cos (2 \pi x) \sin (2 \pi y). 
\]
The parameters $\delta =0.025, \varepsilon =0.3$ and $\gamma=2$ give rise to the initial conditions which are depicted in Figure \ref{fig:Ex4_1}. Similar parameters can be found in \cite[Section 4.1]{HLRZ06}.

\captionsetup[subfigure]{labelformat=empty}
\begin{figure}[!htb]
\centering
\subfloat[]{
\includegraphics[scale=0.3]{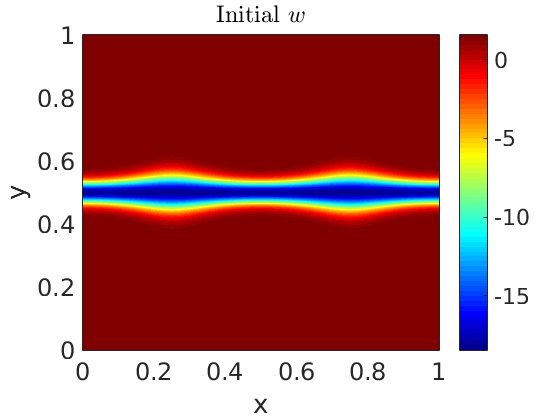}} 
\subfloat[]{
\includegraphics[scale=0.3]{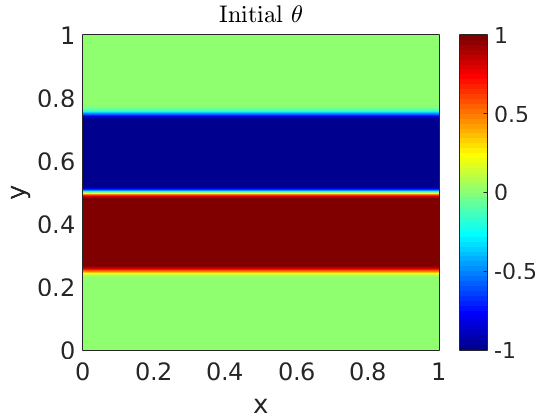}}  \\
\caption{Initial conditions for vorticity $w$ and temperature $\theta$.}
\label{fig:Ex4_1}
\end{figure} 

We apply the DgPC with the following parameters: $K=4$, $N=2$, $S = 2 \times 10^5$, $T=1$, $\Delta t =0.1$, $M=2^7$ and $D= 4,6,8$. Sparse indices $r$ are chosen as: 
\begin{enumerate}[i)]

\item if $D=4$, then $r=(2,1,2,1,2,2,1,1)$,
and if $|\alpha|=2$, we set $(2,\cdot,2,\cdot,2,2)$.  This results in 19 coefficients in the PCE. 

\item if $D=6$, then $r=(2,1,2,1,2,2,1,1,1,1)$,
and if $|\alpha|=2$, we set $(2,\cdot,2,\cdot,2,2,1,1)$.  This results in 30 coefficients in the PCE. 

\item if $D=8$, then $r=(2,1,2,1,2,2,2,2,2,1,1,1)$,
and if $|\alpha|=2$, we set $(2,\cdot,2,\cdot,2,2,2,2,2)$ .   This results in 41 coefficients in the PCE.
\end{enumerate}
We note that sparser sets of indices can be chosen more aggressively in applications to provide faster offline and online computations. 
  
The first three moments of vorticity and temperature are plotted in Figure \ref{fig:Ex4_2}. Higher order moments are centered. Roll-up of the fluid is clearly observed in the mean temperature. Due to the Kelvin-Helmholtz instability and the structure of the initial vorticity, the thin shear layer evolves, rolls up and eventually forms two vortices concentrated at the same locations of sinusoidal perturbations; see \cite{LSM93,HLRZ06,Luo} for the previous results and discussions. 
 
\begin{figure}[!htb]
\centering
\subfloat[]{
\includegraphics[scale=0.3]{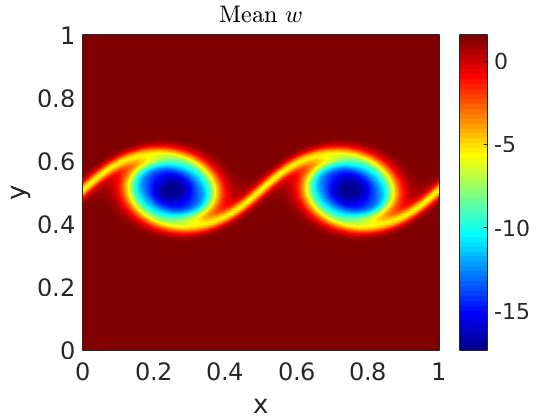}}  
\subfloat[]{
\includegraphics[scale=0.3]{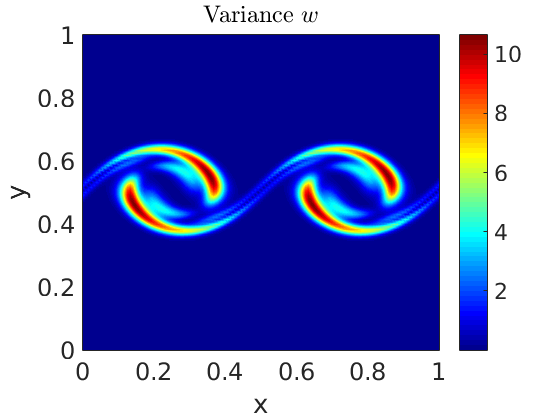}}  
\subfloat[]{
\includegraphics[scale=0.3]{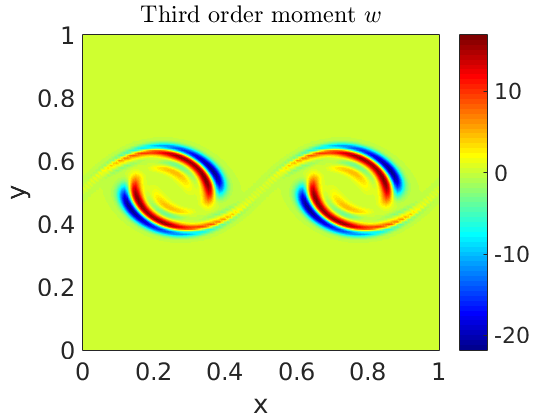}}  \\
\subfloat[]{
\includegraphics[scale=0.3]{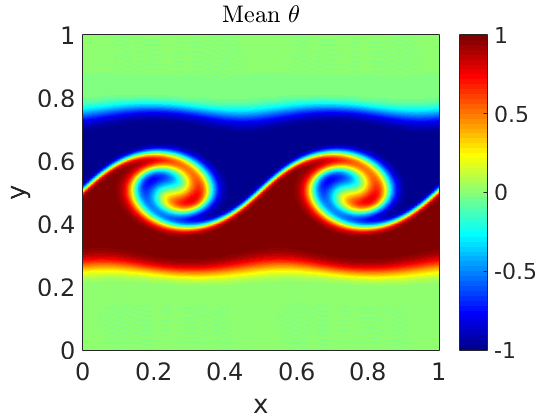}} 
\subfloat[]{
\includegraphics[scale=0.3]{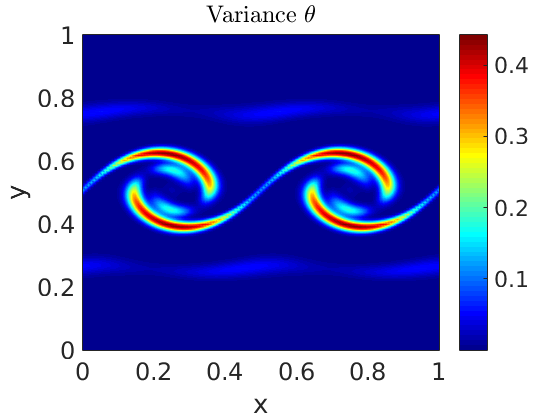}} 
\subfloat[]{
\includegraphics[scale=0.3]{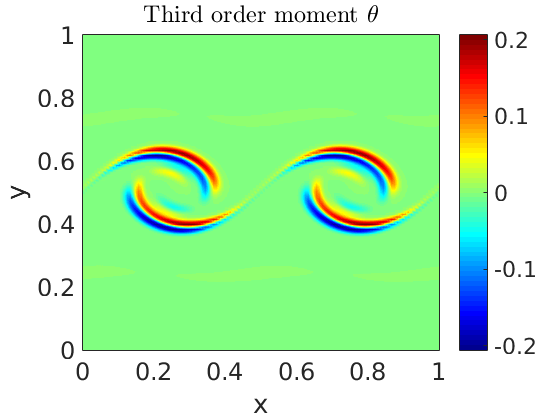}}   
\caption{Moments of vorticity $w$ and temperature $\theta$ obtained by DgPC with $D=8$ at $T=1$.}
\label{fig:Ex4_2}
\end{figure}

To make comparisons as meaningful as possible, we use second order integration methods, namely weak Runge-Kutta and predictor corrector, with the same time-step $dt=0.001$ in both Monte Carlo and DgPC algorithms. Diffusion terms are integrated analytically by an exponential integrator scheme. The number of samples used in MC are $N_{\mbox{samp}}=1000,5000,10000$, and the corresponding algorithms will be denoted by MC1, MC2 and MC3, respectively. Relative $L^2$ errors and computational times are computed using the algorithm MC3 as the ``exact" solution.  

Relative  errors for the moments of the vorticity $w$ are given in Table \ref{table:rel_err_vor}. In this implementation, the algorithm assembles the full covariance matrix at each restart using \eqref{eq:cov_pce} and uses the Arnoldi method to compute the largest $D$  eigenvalues. We found that error levels and convergence behavior for the moments of the temperature $\theta$  were very similar and hence are not displayed. It can be observed that using $D=8$ random modes in DgPC with second degree polynomials yields a similar accuracy to that achieved in MC2. Convergence order of DgPC in terms of the parameter $D$ seems to be at least quadratic whereas convergence of MC is approximately of $O(\sqrt{N_{\mbox{samp}}})$; especially for higher order moments. 

\begin{table}[!htb]
\centering
\begin{tabular}{r|c|c|c |c|c}
Algorithm&          Mean & Variance & $3${rd} order & $4${th} order & Time ratio  \\ \hline
DgPC: $D=4$ & 8.2E-3  &1.58E-1  &7.50E-1 & 4.90E-1 & 0.025 \\ 
DgPC: $D=6$ &  2.7E-3 &3.57E-2 & 1.55E-1& 8.34E-2& 0.041\\ 
DgPC: $D=8$ & 2.1E-3 & 2.64E-2  & 7.01E-2  & 5.98E-2&0.073 \\ 
MC1 & 4.3E-3 & 7.26E-2 &1.42E-1 & 1.16E-1 &0.1\\
MC2 & 3.4E-3  & 2.90E-2 & 5.40E-2  & 4.65E-2 & 0.5 \\ \hline
\end{tabular} 
\caption{Relative errors for moments of the vorticity $w$ at $T=1$ and timing. Exact solution is taken as algorithm MC3.}
\label{table:rel_err_vor}
\end{table}

For comparison, we now apply the random matrix technique discussed in Section \ref{sec:KLexpansion} to compute the KL expansion at each restart. Recall that this technique does not require assembling any covariance matrix and is therefore memory efficient. Table \ref{table:rel_err_fourier_truncation} shows the resulting errors of DgPC algorithm using this technique with the target rank parameter $l=D+p$, where we used $p=10$ for the oversampling parameter. Comparing Table \ref{table:rel_err_vor} and  \ref{table:rel_err_fourier_truncation}, we observe that the error levels remain comparable while the computational costs are not. Elapsed times are approximately divided by $8,4$ and $2$ for the degrees of freedom $D=4,6$ and $D=8$, respectively. For $D=8$, this shows that the (total) algorithm now runs in about half the time. Specifically for the KL step, we  found computational times are reduced by approximately $1000$. This simulation confirms that the computation of the KL expansion becomes a serious computational bottleneck in high dimensions when the covariance matrices are assembled (accounting for half of the time of the full algorithm).  Thanks to the random projection method, the KL expansion no longer constitutes a computational bottleneck for the SNS simulations presented here. The computational costs are now mostly dominated by time evolution and restart procedures.

\begin{table}[!htb]
\centering
\begin{tabular}{r|c|c|c |c|c}
Algorithm&          Mean & Variance & $3${rd} order & $4${th} order & Time ratio  \\ \hline
DgPC: $D=4$ & 8.2E-3  &1.59E-1  & 7.57E-1 & 4.98E-1 & 0.003  \\
DgPC: $D=6$ &  2.8E-3  &3.73E-2 & 1.58E-1& 8.25E-2& 0.010\\  
DgPC: $D=8$  & 2.1E-3 &2.58E-2  & 6.67E-2  & 5.90E-2& 0.035  \\ 
\hline
\end{tabular} 
\caption{Elapsed times and relative errors for moments of vorticity $w$ at $T=1$. The random projection technique (Section \ref{sec:KLexpansion}) with the parameter $l=D+10$ is used to accelerate computation of the KL expansion.}
\label{table:rel_err_fourier_truncation}
\end{table}

\end{exmp}

\begin{exmp} \label{ex:randomvis} \rm 

Using the same scenario as in the previous example, we consider a stochastic viscosity $\nu = U(0.0002,0.0004)$ and set $\mu=\nu$. Since the Monte Carlo simulation takes a very large amount of time to compute, we restrict ourselves to the final time $T=0.5$ and the mesh size $M=2^6$. Monte Carlo algorithms MC1, MC2 and MC3, are executed with the number of samples $100 \times 100$, $200 \times 200$ and $300 \times 300$, respectively. These samples correspond to the samples of Brownian motion and the viscosity. The parameters of DgPC remain the same except we increase $S$ to $3 \times 10^5$ as there is an additional random coefficient in the system. 

\begin{table}[!htb]
\centering
\begin{tabular}{r|c|c|c |c|c}
Algorithm&          Mean  & Variance & $3${rd} order & $4${th} order & Time ratio  \\ \hline
DgPC: $D=4$ & 3.26E-4  & 2.42E-2  & 3.05E-1 & 5.77E-2  & 0.0009 \\ 
DgPC: $D=6$ & 2.85E-4 & 1.63E-2  & 1.95E-1  &  4.89E-2& 0.0025\\ 
DgPC: $D=8$ & 2.74E-4 & 4.45E-3  &  6.27E-2 & 3.30E-2  & 0.0067\\ 
MC1 & 2.60E-3& 2.29E-2  & 9.59E-2 & 5.25E-2  & 0.11\\
MC2 & 1.11E-3 & 8.51E-3  &4.10E-2  & 2.07E-2 & 0.44  \\
\hline
\end{tabular}
\caption{ Relative errors for moments of vorticity $w$ at $T=0.5$. Elapsed times are compared to Algorithm MC3.  }
\label{table:rel_err_vor_rand_vis}
\end{table}

Table \ref{table:rel_err_vor_rand_vis} exhibits the relative errors of DgPC for the vorticity using the random matrix approach for the KL expansion. Comparing Table \ref{table:rel_err_vor_rand_vis}  and \ref{table:rel_err_fourier_truncation}, we see that relative elapsed times of DgPC with respect to MC3 are further improved. Additional randomness for MC means an extra dimension to sample from, whereas for DgPC, it means an extra variable that needs to be compressed into the modes $\bta$. Since the dynamics crucially depend on the behavior of the viscosity, using few realizations for viscosity sampling in MC is not recommended. Iin this setting, we found that our MC simulations demanded a high CPU time compared to DgPC.   Note, however, that viscosity sampling could clearly be performed in parallel in a MC framework---something that is not as easily feasible in the PCE setting.

\end{exmp}

\begin{exmp} \rm 

The preceding simulations were concerned with short time evolutions of SNS and comparisons of the proposed algorithm with a Monte Carlo method. Numerical results for reasonably short time computations indicated that our algorithm achieved a similar accuracy compared to MC typically for a smaller computational cost.

We are now interested in long time simulations and convergence to steady states. Since there is no random forcing acting upon the temperature equation in \eqref{eq:sns_vor}, the (uncoupled) temperature diffuses to zero quickly. Therefore, we only solve the vorticity equation in the system \eqref{eq:sns_vor}. 

The following numerical experiment considers the vorticity equation with a deterministic viscosity $\nu=0.00055$ and a spatial forcing as described in Example \ref{ex:sns_short1}. The parameters of the simulation are  $M=2^6$, $S=3 \times 10^5$, $T=288$ and $\Delta t= 0.12$. PC expansions with thirty number of terms are employed on each subinterval.  The four-step Adams predictor-corrector method is used for the time integration.   

Figure \ref{fig:Ex7_1}  shows three different initial conditions for the vorticity. The first layer is supported around $x=0.5$ while the others are aligned horizontally. Widths of all layers are widened and different sinusoidal perturbations are considered. 

\captionsetup[subfigure]{labelformat=empty}
\begin{figure}[!htb]
\centering
\subfloat[]{
\includegraphics[scale=0.3]{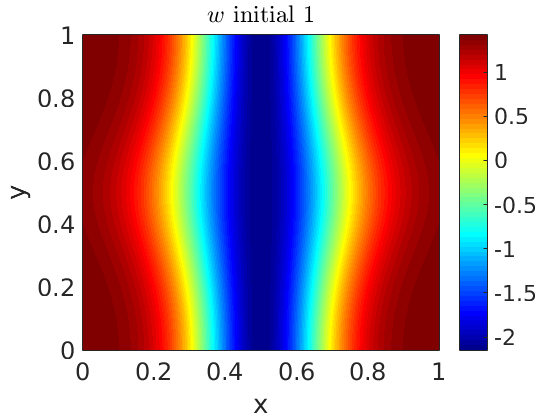}} 
\subfloat[]{
\includegraphics[scale=0.3]{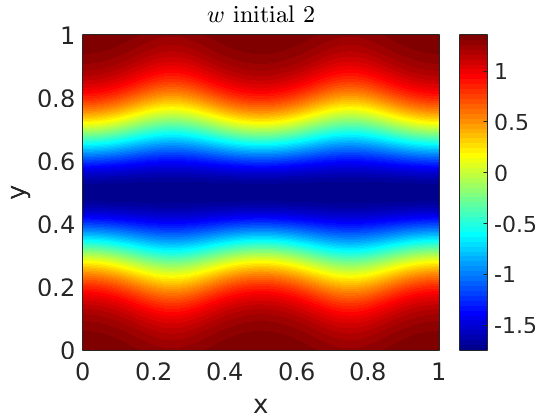}} 
\subfloat[]{
\includegraphics[scale=0.3]{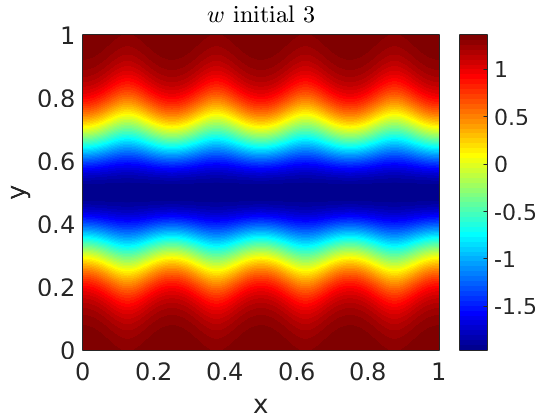}} 
\caption{Different initial conditions for the vorticity $w$.}
\label{fig:Ex7_1}
\end{figure} 

In Figure \ref{fig:Ex7_2}, we show the $L^2$-norm of the successive differences of the first two moments in time. Each column represents one of the initial conditions presented in the corresponding column in Figure \ref{fig:Ex7_1}. After a (very) long time, the norms of the successive differences drop below $O(10^{-3})$, which (numerically) indicates that statistical moments no longer change significantly in time. In all cases, we found that the dynamics converged to the same state, which is shown in Figure \ref{fig:Ex7_3}. Notice that the invariant measure is a non-Gaussian random field and the moments have oscillations in the $x$ variable. We also see that high variance regions correspond to where the mean fields display peaks.  Based on these findings for this scenario, we assert that the dynamics converge to an invariant measure which is numerically captured in the long-time by the DgPC algorithm.

\captionsetup[subfigure]{labelformat=empty}
\begin{figure}[!htb]
\centering 
\subfloat[]{
\includegraphics[scale=0.3]{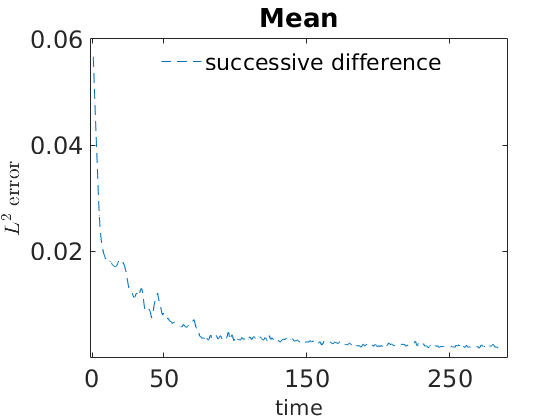}} 
\subfloat[]{
\includegraphics[scale=0.3]{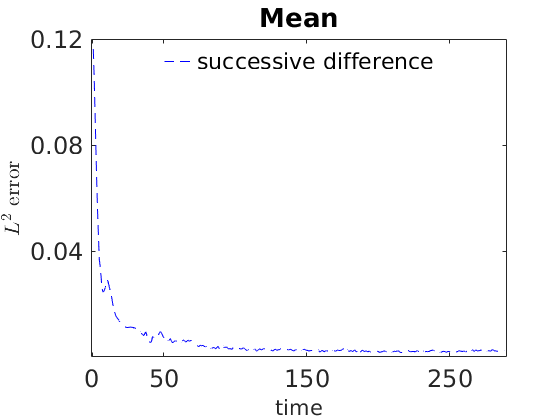}} 
\subfloat[]{
\includegraphics[scale=0.3]{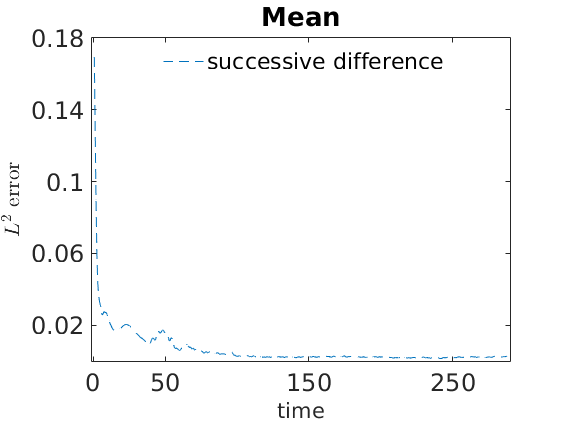}}  \\
\vspace{-0.3cm}
\subfloat[]{
\includegraphics[scale=0.3]{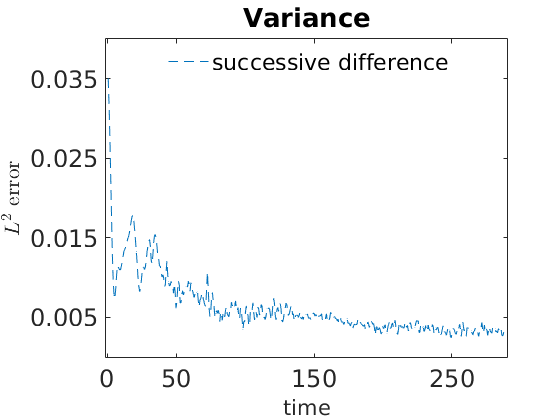}} 
\subfloat[]{
\includegraphics[scale=0.3]{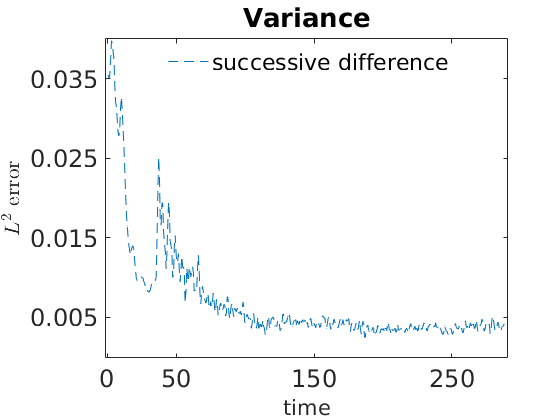}} 
\subfloat[]{
\includegraphics[scale=0.3]{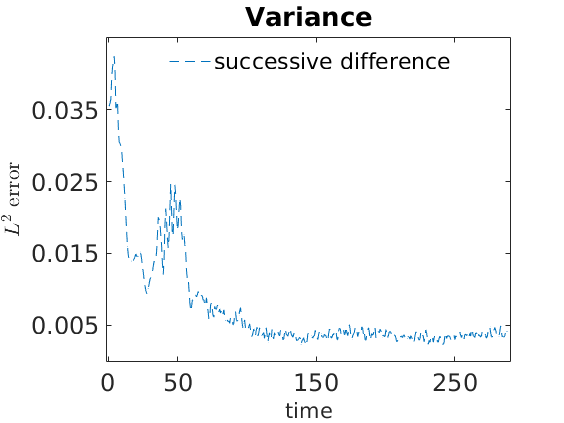}}
\caption{$L^2$-norm of successive differences of moments using three different initial conditions. Each column corresponds to an initial condition depicted in Figure \ref{fig:Ex7_1}. }
\label{fig:Ex7_2}
\end{figure}

 \captionsetup[subfigure]{labelformat=empty}
\begin{figure}[!htb]
\centering
\subfloat[]{
\includegraphics[scale=0.3]{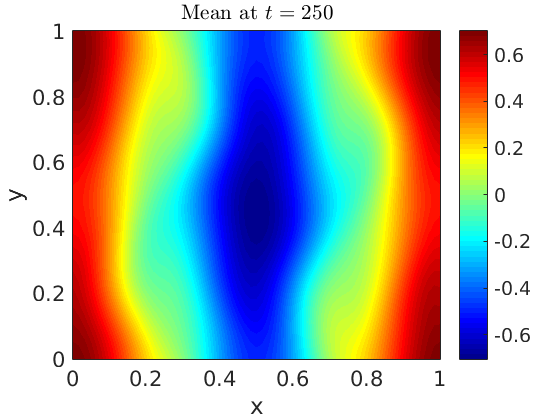}} 
\subfloat[]{
\includegraphics[scale=0.3]{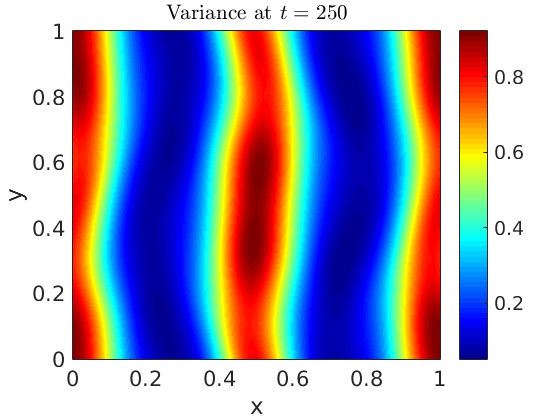}}
\subfloat[]{
\includegraphics[scale=0.3]{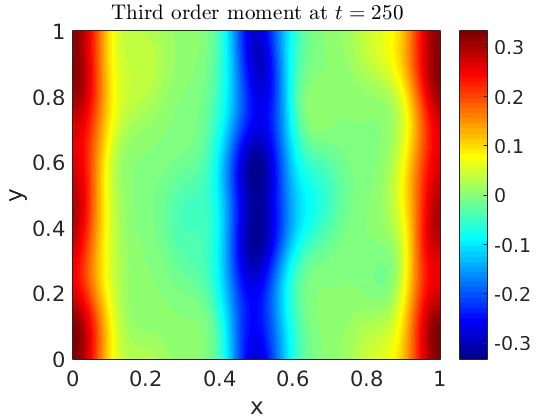}}  
\caption{Statistical moments of the invariant measure of the vorticity at time $t=250$ obtained by DgPC with three different initial conditions depicted in Figure \ref{fig:Ex7_1}.}
\label{fig:Ex7_3}
\end{figure}

\end{exmp}

\begin{remark} \rm
 In long time computations, the MC method  usually requires the propagation of many realizations in time, which renders the method hardly affordable in some cases. However, if the dynamical system possesses a unique ergodic invariant measure, the MC method may be carried out to sample such a measure by considering a single, very long time MC realization, which repeatedly visits the whole state space. While carrying out such a sampling is also computationally expensive, it is likely to compare favorably to our DgPC algorithm in this case.
 
In general, ergodicity or uniqueness of an invariant measure may not be known or not hold for complicated physical dynamics (e.g., invariant measures parametrized by a random parameter as in Example \ref{ex:randomvis}). In such cases, our algorithm offers a viable alternative to the MC method to capture the long-term dynamics by providing statistical information resulting from the expansion coefficients.

\end{remark}

\section{Summary} \label{sec:summary}

We have presented a PC-based algorithm, called Dynamical generalized Polynomial Chaos, to tackle long-time integration and high dimensional randomness in the context of PDEs with Markovian forcing. To deal with these challenges, DgPC uses a restart procedure, which constructs PCEs dynamically based on the polynomials of the projections of the solution, the random forcing and the random parameters found in the equation. The Karhunen--Loeve expansion is employed at each restart to find a representation of the solutions which correspond to low-dimensional dynamics of the underlying physical processes. The relevant modes are then incorporated in a PCE for the future evolution. Using sparse multi-index sets and frequent restarts, the algorithm provides an efficient way to capture the solutions in a fairly sparse random bases in terms of orthogonal polynomials of dynamically evolving measures.

The main computational bottlenecks of the algorithm are the simulation of the deterministic evolution equation, the KL expansion, and the computation of moments. The cost of the deterministic evolution is dictated by the complicated nature of the SPDEs. The KL expansion is an expensive dimensionality reduction technique. Since the algorithm constructs PCEs online, KL expansions (or other dimensionality reductions methods) are unavoidable at each restart. We found that for large covariance matrices, the KL cost was drastically reduced when the covariance matrix was estimated by a low-rank approximation obtained by random projections. The estimation of the orthogonal polynomials and corresponding triple products of evolving arbitrary measures is also a costly step as in most PC-based methods.

Using a 1D randomly forced Burgers equation and a 2D stochastic Navier--Stokes system, we provided several numerical simulations for both short- and long-time solutions. We compared the accuracy and computational time of the algorithm to the standard Monte Carlo method and found that the proposed algorithm achieved similar error levels for a (generally significant) lower computational cost in most cases. The substantial speed-up of DgPC is especially promising when the equation contains additional, time-independent, random contributions, which is one of the main reasons to use PC--based methods in general. To demonstrate the efficiency of the algorithm for long time simulations, we computed invariant measures for both equations, which is not a trivial task for two dimensional Navier--Stokes systems.

Other methods such as the recent  Multilevel Monte Carlo techniques offer improvements over the standard MC methods \cite{G08}. The restarting step of our algorithm remains expensive computationally, especially for problems in two (or three) spatial dimensions. However, the restart method provides a viable means to keep the number of random variables to reasonable levels. Its ability to compute statistical properties of long-time evolutions for fairly complicated equations is quite promising.

\section*{Acknowledgement}
The authors would like to thank the reviewers for their critical reading of the manuscript, and for multiple remarks and suggestions that improved the presentation of our method.
This work was partially funded by NSF Grant DMS-1408867 and ONR Grant N00014-15-1-2679.

\bibliographystyle{plain}

\bibliographystyle{plain}

\end{document}